\documentclass[a4paper,12pt]{article}

\usepackage{mymacros}

\newcommand{\CY}{\mathrm{CY}}
\newcommand{\GLt}{\widetilde{\GL}}

\newcommand{\res}{\mathrm{res}}
\newcommand{\VC}{\operatorname{VC}}

\newcommand{\Gmax}{{G_{\mathrm{max}}}}
\newcommand{\Gbarmax}{{\Gbar_{\mathrm{max}}}}
\newcommand{\Gvmax}{{\Gv_{\mathrm{max}}}}

\newcommand{\Eff}{\operatorname{Eff}}
\newcommand{\Euler}{\operatorname{Euler}}
\newcommand{\td}{\operatorname{td}}
\newcommand{\ch}{\operatorname{ch}}

\newcommand{\Pictor}{\Pic_{\mathrm{tor}}}
\newcommand{\Piccor}{\Pic_{\mathrm{cor}}}

\newcommand{\cDv}{{\Check{\cD}}}
\newcommand{\cDbar}{{\overline{\cD}}}

\newcommand{\HA}{\scrH_A}
\newcommand{\HB}{\scrH_B}
\newcommand{\HAZ}{H_{A, \, \bZ}}
\newcommand{\HBZ}{H_{B, \, \bZ}}

\newcommand{\HBQ}{H_{B, \, \bQ}}

\newcommand{\HAC}{H_{A, \, \bC}}

\newcommand{\Mir}{\operatorname{Mir}}
\newcommand{\amb}{{\mathrm{amb}}}
\newcommand{\reg}{{\mathrm{reg}}}
\newcommand{\vc}{{\mathrm{vc}}}

\newcommand{\NS}{\operatorname{NS}}
\newcommand{\Stab}{\operatorname{Stab}}

\newcommand{\PD}{\operatorname{PD}}

\title{Mirror symmetry and K3 surfaces}
\author{Kazushi Ueda}
\date{}
\begin{document}

\maketitle

\begin{abstract}
We review some of the interplay
between mirror symmetry and K3 surfaces.
\end{abstract}


\section{Introduction}

\emph{Mirror symmetry} is a mysterious relationship,
originally suggested by string theorists,
between symplectic geometry of one Calabi-Yau manifold
and complex geometry of another Calabi-Yau manifold.
One of the earliest prediction in mirror symmetry
is that Calabi-Yau 3-fold should come in pairs
$(Y, \Yv)$
in such a way that
their Hodge numbers satisfy
$h^{1,1}(Y) = h^{1,2}(\Yv)$ and
$h^{1,2}(Y) = h^{1,2}(\Yv)$.
Note that $h^{1,1}(Y)$ is the number of parameters
for the K\"{a}hler structures, and
$h^{1,2}(Y)$ is the dimension
of the moduli space of complex structures.
More generally,
a pair $(Y, \Yv)$ of Calabi-Yau $n$-folds
is said to be a \emph{topological mirror pair}
if their Hodge numbers satisfy
\begin{align} \label{eq:tms}
 h^{i,j}(Y) = h^{i,n-j}(\Yv)
\end{align}
for any $i, j$.
Although topological mirror symmetry is much weaker
than other versions of mirror symmetry,
construction of a topological mirror partner
of a given Calabi-Yau manifolds is already a highly non-trivial problem.
One subtlety is the existence of rigid Calabi-Yau manifolds,
as the mirror partner of such manifolds
can not exist as a Calabi-Yau manifold.
Another subtlety comes in
if one allows Calabi-Yau $n$-folds to have singularities,
which motivated the theory of stringy Hodge numbers and
motivic integration.

For K3 surfaces,
topological mirror symmetry seems to be trivial
at first sight,
since every K3 surface has the identical Hodge numbers.
The moduli space of complex structures
and the moduli space of K\"{a}hler structures are somehow `mixed',
as they both live in $H^2(Y; \bC)$.
In a sense, it is more natural to work with
the moduli space of hyperK\"{a}hler structures,
instead of those of complex structures
and K\"{a}hler structures.
However,
`topological' mirror symmetry for K3 surfaces
can be formulated,
not as an exchange of the Hodge numbers,
but as an exchange of the algebraic lattice
and the transcendental lattice
inside the total cohomology group.
As such, it is much more subtle
than topological mirror symmetry
for Calabi-Yau 3-folds,
partly because lattices have more structures
than Hodge numbers, and
partly because these lattices are sensitive
to the complex structure of the K3 surface.

\emph{Classical mirror symmetry} is
the mysterious relationship
between Gromov-Witten invariants of $Y$ and
period integrals of its mirror manifold $\Yv$.
It is also known as \emph{Hodge-theoretic mirror symmetry},
since it can be formulated
as an isomorphism
between two variations of Hodge structures.
One of them,
called the A-model VHS,
is defined
on the `moduli space of complexified K\"{a}hler structures' of $Y$,
and encode the information of Gromov-Witten invariants of $Y$.
There is no satisfactory definition
of the moduli space of complexified K\"{a}hler structures of $Y$;
the space of stability conditions
on the derived category of coherent sheaves
\cite{Bridgeland_SCTC, Bridgeland_SCK3}
is expected to be the universal cover of this moduli space.
The other one,
called the B-model VHS,
is the variation of Hodge structures
on the moduli space of complex structures on $\Yv$
defined by the usual Hodge theory.
Classical mirror symmetry started
with the prediction
that the numbers of rational curves
on a general quintic hypersurface in $\bP^3$
can be computed
by period integrals of its mirror family
\cite{CdGP},
and attracted much attention
from mathematicians.

Deformation invariance of Gromov-Witten invariants implies
that Gromov-Witten invariants of K3 surfaces are trivial,
since a generic K3 surface does not have
any holomorphic curve at all.
Hence classical mirror symmetry for K3 surfaces
reduces essentially to the study of period maps.
In particular,
the Yukawa coupling can be identified
with the cup product on the mirror
\cite{Dolgachev_MSK3}.

\emph{Homological mirror symmetry}
is introduced by Kontsevich
\cite{Kontsevich_HAMS}
to give a deeper understanding of mirror symmetry.
A pair $(Y, \Yv)$
of Calabi-Yau manifolds
is a \emph{homological mirror pair}
if one has an equivalence
\begin{align} \label{eq:hms}
 D^\pi \Fuk Y \cong D^b \coh \Yv
\end{align}
of derived categories.
It is hard to find a homological mirror pair,
and the only known homological mirror pair
of K3 surfaces is the case
when $Y$ is a quartic hypersurface in $\bP^3$,
due to Seidel \cite{Seidel_K3}.
There is a fascinating interplay
between homological mirror symmetry
and monodromy of period maps,
also pioneered by Kontsevich.
It is closely related to stability conditions
on triangulated categories
introduced by Bridgeland
\cite{Bridgeland_SCTC}.

A geometric picture for mirror symmetry
is provided by the \emph{Strominger-Yau-Zaslow conjecture}
\cite{Strominger-Yau-Zaslow},
which states that
any Calabi-Yau manifold has a structure
of a special Lagrangian torus fibration
$\pi \colon Y \to B$,
and its mirror is obtained as
the dual special Lagrangian torus fibration
$\piv \colon \Yv \to B$.
This motivated the construction of mirror manifolds
using integral affine manifolds with singularities
\cite{MR2181810,Gross-Siebert_RAGCG}.
We will not review this here,
and refer the interested reader
to an excellent review
\cite{1212.4220}
and references therein.

This review is organized as follows:
In \pref{sc:sd},
we recall strange duality
and its generalizations,
which are now understood
as incarnations of mirror symmetry.
In \pref{sc:BH},
we discuss transposition mirror construction
by Berglund and H\"{u}bsch.
%
%
In \pref{sc:am},
we discuss mirror symmetry for K3 surfaces
following Aspinwall and Morrison.
%
In \pref{sc:Dolgachev},
we discuss the notion of lattice-polarized K3 surfaces,
which is introduced
by Nikulin
and used by Dolgachev
to study mirror symmetry for K3 surfaces.
In \pref{sc:Batyrev},
we discuss mirror construction
due to Batyrev
using polar duality of reflexive polytopes.
In \pref{sc:classical},
we discuss classical mirror symmetry
for anti-canonical hypersurfaces in toric weak Fano 3-folds.
%
In \pref{sc:D_conj},
we discuss a conjecture of Dolgachev
on the relation between Batyrev mirrors
and Dolgachev mirrors of K3 surfaces.
In \pref{sc:stab},
we discuss Bridgeland stability conditions on K3 surfaces.
%
In \pref{sc:BV},
we discuss mirror construction,
due to Borcea
and Voisin,
for Calabi-Yau 3-folds
associated with 2-elementary K3 surfaces.


\section{Strange duality}
 \label{sc:sd}


\subsection{The modality of a singularity}

Two germs
$(f^{-1}(0), 0)$
and $(g^{-1}(0), 0)$
of hypersurface singularities
defined by convergent power series
$
 f \colon (\bC^n, 0) \to (\bC, 0)
$
and
$
 g \colon (\bC^n, 0) \to (\bC, 0)
$
are \emph{right equivalent}
if there is a holomorphic change of coordinates
$
 \varphi \colon (\bC^n, 0) \to (\bC^n, 0)
$
such that $f = g \circ \varphi$.
If the critical point
of the germ $f \colon (\bC^n, 0) \to (\bC, 0)$
of a holomorphic function is isolated,
then it is right equivalent
to a polynomial of degree at most $\mu+1$
\cite{MR0240826}.
Here $\mu$ is the \emph{Milnor number} of $f$,
defined as the dimension of the Jacobi ring
$
 \bC \{ \! \{ z_1, \ldots, z_n \} \! \} / (\partial_{z_1} f, \ldots, \partial_{z_n} f)
$
of $f$.
This \emph{finite-determinacy} of isolated singularities
allows one to reduce the classification
of isolated critical points of holomorphic functions
up to right equivalence
to a finite-dimensional problem.

Let $G$ be a Lie group acting on a manifold $M$.
The \emph{modality} of a point $f \in M$
is the smallest integer $m$ such that
a sufficiently small neighborhood of $f$
can be covered by a finite number of $m$-parameter families
of orbits.
The modality of the germ of a function
is the modality of its jet
in the space of $k$-jets
with respect to the right action
of the group of coordinate transformations
for sufficiently large $k$.
If $f$ is the defining equation of an isolated singularity,
then the modality of $f$ is called the modality of the singularity.
In other words,
the modality of a singularity
is the minimal number of continuous parameters
needed to parametrize right equivalence classes of singularities
that are close to $p$.
Singularities of modality zero are called \emph{simple singularities},
and have well-known ADE classification.
Next in line come \emph{unimodal singularities},
which are classified by Arnold
\cite{Arnold_ICM}
into an infinite series
\begin{align} \label{eq:Tpqr}
 T_{p,q,r} : x^p + y^q + z^r + a x y z
\end{align}
and 14 exceptional cases.
\begin{table}[t]
\[
\begin{array}{cccccc}
  \toprule
 \text{name} & \text{normal form} & \text{weight} &
 \bsdelta & \bsgamma & \text{dual} \\
  \midrule
 E_{12} & x^3 + y^7 + z^2 + a x y^5
  & (6,14,21;42) & (2,3,7) & (2,3,7) & E_{12} \\
 E_{13} & x^3 + x y^5 + z^2 + a y^8
  & (4,10,15;30) & (2,4,5) & (2,3,8) & Z_{11} \\
 E_{14} & x^3 + y^8 + z^2 + a x y^6
  & (3,8,12;24) & (3,3,4) & (2,3,9) & Q_{10} \\
 Z_{11} & x^3 y + y^5 + z^2 + a x y^4
  & (6,8,15;30) & (2,3,8) & (2,4,5) & E_{13} \\
 Z_{12} & x^3 y + x y^4 + z^2 + a x^2 y^3
  & (4,6,11;22) & (2,4,6) & (2,4,6) & Z_{12} \\
 Z_{13} & x^3 y + y^6 + z^2 + a x y^5
  & (3,5,9,18) & (3,3,5) & (2,4,7) & Q_{11} \\
 W_{12} & x^4 + y^5 + z^2 + a x^2 y^3
  & (4,5,10;20) & (2,5,5) & (2,5,5) & W_{12} \\
 W_{13} & x^4 + x y^4 + z^2 + a y^6
  & (3,4,8,16) & (3,4,4) & (2,5,6) & S_{11} \\
 Q_{10} & x^3 + y^4 + y z^2 + a x y^3
  & (6,8,9;24) & (2,3,9) & (3,3,4) & E_{14} \\
 Q_{11} & x^3 + y^2 z + x z^3 + a z^5
  & (4,6,7;18) & (2,4,7) & (3,3,5) & Z_{13} \\
 Q_{12} & x^3 + y^5 + y z^2 + a x y^4
  & (3,5,6;15) & (3,3,6) & (3,3,6) & Q_{12} \\
 S_{11} & x^4 + y^2 z + x z^2 + a x^3 z
  & (4,5,6,16) & (2,5,6) & (3,4,4) & W_{13} \\
 S_{12} & x^3 y + y^2 z + x z^3 + a z^5
  & (3,4,5;13) & (3,4,5) & (3,4,5) & S_{12} \\
 U_{12} & x^3 + y^3 + z^4 + a x y z^2
  & (3,4,4;12) & (4,4,4) & (4,4,4) & U_{12} \\
  \bottomrule
\end{array}
\]
\caption{14 exceptional unimodal singularities}
\label{tb:EUS}
\end{table}
Here $(p,q,r)$ is a triple of natural numbers
satisfying $1/p+1/q+1/r \le 1$,
and $a$ is a general complex parameter.
The $T_{p,q,r}$-singularity is called
a \emph{simple elliptic singularity}
if $1/p+1/q+1/r = 1$,
and a \emph{cusp singularity}
if $1/p+1/q+1/r < 1$.
\emph{Exceptional unimodal singularities} are listed
in \pref{tb:EUS},
where $a$ is again a general complex parameter.
The defining polynomial is weighted homogeneous
for $a = 0$,
and the weights in \pref{tb:EUS} are the primitive weights
of the variables and the defining polynomial.

\begin{figure}
\centering
\phantom{a} \ \\[-3mm]
\input{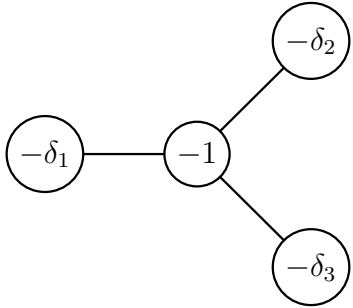}
\phantom{a} \ \\[3mm]
\caption{The diagram $T(\delta_1, \delta_2, \delta_3)$}
\label{fg:dolg_num}
\end{figure}

There are two natural ways to study singularities;
one is to \emph{resolve} the singularity,
and the other is to \emph{deform} the singularity.
%
A resolution of a surface singularity is \emph{good}
if the exceptional locus is a simple normal crossing divisor.
A good resolution is \emph{minimal}
if any contraction of an exceptional curve
gives a non-good resolution.
The minimal good resolution of an exceptional unimodal singularity
consists of four rational curves;
one is an exceptional curve of the first kind,
and the others are mutually disjoint rational curves,
each intersecting the first curve in one point.
\pref{fg:dolg_num} shows the dual graph
of the exceptional divisor.
The self-intersection numbers
$\bsdelta = (\delta_1, \delta_2, \delta_3)$
of the three exceptional curves
is called the \emph{Dolgachev number} of the singularity.

\begin{figure}
\centering
\input{ext_dynkin_diagram.pst}
\caption{The diagram $\That(\gamma_1, \gamma_2, \gamma_3)$}
\label{fg:Dynkin-Milnor}
\end{figure}

The \emph{Milnor fiber} of a singularity
is the intersection
$
 f^{-1}(\epsilon) \cap B_\delta
$
of the deformation $f^{-1}(\epsilon)$
of the singularity
with a sufficiently small ball
$
 B_\delta = \{ (x,y,z) \in \bC^3 \mid |x|^2+|y|^2+|z|^2 \le \delta \}.
$
Here $\epsilon$ is a sufficiently small number
which may depend on $\delta$.
The diffeomorphism type of the Milnor fiber
is independent of the choice of $\delta$ and $\epsilon$.
The Milnor fiber is homotopy-equivalent
to the bouquet $(S^2)^{\vee \mu}$ of $\mu$ 2-spheres,
where $\mu$ is the Milnor number
\cite{MR0239612}.
The middle-dimensional homology group of the Milnor fiber,
equipped with the intersection form,
is called the \emph{Milnor lattice}.
The Milnor lattice has a distinguished basis
$\{ \alpha_i \}_{i=1}^\mu$
consisting of \emph{vanishing cycles},
which is well-defined up to an action of a braid group.
By a suitable choice of a distinguished basis of vanishing cycles,
the Coxeter-Dynkin diagram of an exceptional unimodal singularity
can be written in the form $\That_\bsgamma$
given in \pref{fg:Dynkin-Milnor}
\cite{MR0367274}.
The triple $\bsgamma = (\gamma_1, \gamma_2, \gamma_3)$
is called the \emph{Gabrielov number} of the singularity.
%
The \emph{strange duality} is an obersevation
by Arnold that
exceptional unimodal singularities come in pairs
in such a way that the Dolgachev number
and the Gabrielov number are interchanged.
For more on foundations of singularity theory,
one can see \cite{MR777682,MR966191,MR0239612}
and references therein.

\subsection{Triangle singularities}

For a sequence
$(p, q, r) \in (\bZ^+)^3$
of positive integers,
the \emph{triangle group}
\begin{align}
 \Delta_{p, q, r}
  := \la a, b, c \mid a^2 = b^2 = c^2 = (ab)^p = (bc)^q = (ca)^r = 1 \ra
\end{align}
is the group generated by reflections
in the triangle of angles $\pi/p$, $\pi/q$ and $\pi/r$.
The triangle group is called {\em spherical}, {\em Euclidean}, or
{\em hyperbolic}
depending on whether
$1/p + 1/q + 1/r$
is greater than, equal to, or less than one respectively.
Let
$
 \Gamma_{p, q, r} \subset \Delta_{p,q,r}
$
be the \emph{von Dyck group},
which is the index two subgroup
of the triangle group
consisting of compositions of even numbers of reflections.
It is written as
\begin{align}
 \Gamma_{p, q, r}
  := \la \xbar, \ybar, \zbar \mid \xbar^p
   = \ybar^q = \zbar^r = \xbar \ybar \zbar = 1 \ra
\end{align}
where the inclusion $\Gamma_{p, q, r} \hookrightarrow \Delta_{p,q,r}$
sends $\xbar$, $\ybar$, and $\zbar$
to $ab$, $bc$ and $ca$ respectively.
The sequence $(p,q,r)$ is called the \emph{signature}
of the triangle group.

\subsubsection{Spherical case}

\begin{table}
\[
\begin{array}{cccccc}
  \toprule
 \text{type} & A_{p+q} & D_{n+2} & E_6 & E_7 & E_8 \\
  \midrule
 \text{signature} &
  (1,p,q) & (2,2,n) & (2,3,3) & (2,3,4) & (2,3,5) \\
  \bottomrule
\end{array}
\]
\caption{Spherical signatures}
\label{tb:spherical_wt}
\end{table}

\begin{figure}[t]
\centering
\input{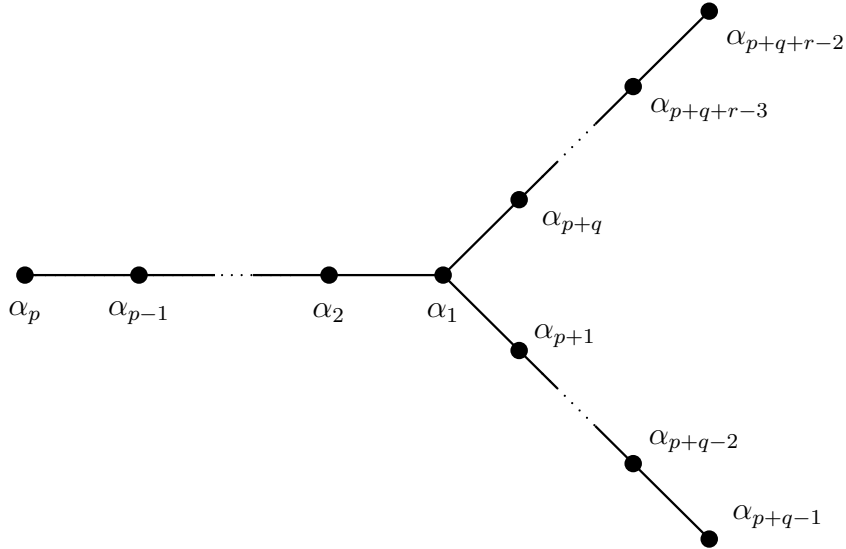}
\caption{The diagram $T_{p,q,r}$}
\label{fg:Tpqr}
\end{figure}

Spherical signatures are classified
in \pref{tb:spherical_wt}.
The von Dyck group $\Gamma_{p,q,r}$ is a polyhedral group,
which is a finite subgroup
of $\PSL_2(\bC)$ (or $SO_3(\bR)$).
The Kleinian group $\Pi_{p,q,r}$
obtained as the pull-back of $\Gamma_{p,q,r}$
by the universal covering map $\SL_2(\bC) \to \PSL_2(\bC)$
is a binary polyhedral group,
which is described as an abstract group as
$$
 \Pi_{p,q,r} = \la x, y, x \mid x^p = y^q = z^r = x y z \ra.
$$
Let
$
 \bX = \bX_{p,q,r} := [\bP^1 / \Gamma_{p,q,r}]
$
be the quotient stack, and
$
 \bK = \bK_{p,q,r}
$
be the total space of the canonical (orbi-)bundle.
The coarse moduli space $X$ of $\bX$ is isomorphic to $\bP^1$,
and $\bX$ has three orbifold points
with stabilizer groups $\bZ / p \bZ$,
$\bZ / q \bZ$ and $\bZ / q \bZ$.
Let
\begin{align*}
 R
  = H^0(\scO_\bK)
  = \bigoplus_{k=0}^\infty H^0(\scO_\bX(-k \omega_\bX)),
\end{align*}
be the anticanonical ring,
which is the ring of regular functions on $\bK$.
The spectrum $\Ybar := \Spec R$ of $R$ gives
the corresponding Kleinian singularity
$
 \Ybar \cong \bC^2 / \Pi_{p,q,r},
$
and the coarse moduli space $K$ of $\bK$
is a partial resolution of $Y$
with three simple singularities
of types $A_{p-1}$, $A_{q-1}$ and $A_{r-1}$.
The minimal resolution $Y$ of $K$
is the minimal resolution of $\Ybar$,
and the exceptional divisor is a tree of $(-2)$-curves
where three chains of $(-2)$-curves
consisting of $p-1$, $q-1$ and $r-1$ rational curves
is connected to the $(-2)$-curve
obtained as the strict transform
of the zero-section in $K$.
The dual graph of the exceptional curves
in the resolution $Y \to \Ybar$ is the $T_{p,q,r}$-graph
shown in \pref{fg:Tpqr}.
One has a derived equivalence
$$
 D^b \coh \bK \cong D^b \coh Y
$$
by Kapranov and Vasserot
\cite{Kapranov-Vasserot}.

%

\subsubsection{Euclidean case}

\begin{table}
\[
\begin{array}{cccc}
  \toprule
 \text{type} & \Etilde_6 & \Etilde_7 & \Etilde_8 \\
  \midrule
 \text{signature} & (3,3,3) & (2,4,4) & (2,3,6) \\
  \text{weight} & (1,1,1) & (1,1,2) & (1,2,3) \\
 E \cdot E & -3 & -2 & -1 \\
  \bottomrule
\end{array}
\]
\caption{Euclidean signatures}
\label{tb:parabolic_wpl}
\end{table}

Euclidean signatures are classified in \pref{tb:parabolic_wpl},
which correspond to the tessellation of $\bC$
by triangles of angles $2 \pi/p$, $2 \pi/q$, and $2 \pi/r$.
The von Dyck group is a subgroup of
$
 SO_2(\bR) \ltimes \bR^2
  \cong U(1) \ltimes \bC
$
acting naturally on $\bC$.
The quotient $E = \bC/\Gamma_{p,q,r}$ is a smooth elliptic curve,
which can be written as
$
 \{ [x:y:z] \in \bP(a,b,c) \mid x^p + y^q + z^r = 0 \}
$
where the weights $(a,b,c)$ are given by
$(1,1,1)$, $(1,1,2)$, and $(1,2,3)$ respectively.

Recall that a normal surface singularity is a {\em simple elliptic singularity}
if the exceptional divisor $E$
of the minimal resolution
is a smooth elliptic curve.
Simple elliptic singularities,
which are also hypersurface singularities,
are classified by Saito \cite{KSaito_EES}
into three types $\Etilde_6$,
$\Etilde_7$ and $\Etilde_8$.
They are $T_{p,q,r}$-singularities
\eqref{eq:Tpqr}
for $(p,q,r) = (3,3,3), (2,4,4)$, and $(2,3,6)$,
which are triangle singularities
if $a=0$.
The triple $(p,q,r)$ is the signature
of the corresponding triangle group.

\subsubsection{Hyperbolic case}

A hyperbolic von Dyck group is a Fuchsian group,
i.e., a discrete subgroup of $\PSL_2(\bR) \cong \Aut \bH$.
A holomorphic function
$
 f \colon \bH \to \bC
$
is an \emph{automorphic form of weight $k$}
with respect to a Fuchsian group $\Gamma$ if
\begin{align}
 f \lb \frac{a z + b}{c z + d} \rb = (c z + d)^k f(z),
  \quad \forall \begin{pmatrix} a & b \\ c & d \end{pmatrix} \in \Gamma,
\end{align}
and $f$ is holomorphic at the cusp.
The orbifold quotient
$
 \bX^\circ = [\bH /\Gamma]
$
can be compactified to a smooth orbifold $\bX$
in such a way that $\bX \setminus \bX^\circ$ has the trivial stabilizer.
The space $A_k(\Gamma)$ of automorphic forms of weight $k$
can be identified with the space
\begin{align}
 H^0(\bX, (T^* \bX)^{\otimes k})
  \cong H^0(\bH, (T^* \bH)^{\otimes k})^\Gamma
\end{align}
of sections of the $k$-th tensor power of the cotangent bundle.
Let
$
 A(\Gamma)
  = \bigoplus_{n=0}^\infty A_n(\Gamma)
$
be the ring of automorphic forms.
The spectrum $\Vbar_\Gamma = \Spec A(\Gamma)$
is a cone over the orbifold $\bX$, and
has an isolated singularity at the origin
called a \emph{hyperbolic triangle singularity}.
It is known by Dolgachev
that a hyperbolic triangle singularity
is a hypersurface singularity
if and only if it is an exceptional unimodal singularity
\cite{MR0568900}.
The signature $(p,q,r)$ of the von Dyck group coincides
with the Dolgachev number of the corresponding
exceptional unimodal singularity.

\subsection{Strange duality and K3 surfaces}
 \label{sc:SD_K3}

Pinkham \cite{Pinkham_strange-duality}
and Dolgachev and Nikulin
\cite{Dolgachev_IQF,Nikulin_ISBF}
gave the following beautiful interpretation of the strange duality
in terms of the Picard lattices and the transcendental lattices
of K3 surfaces
obtained as compactifications of the Milnor fibers.

Let $\Gamma$ be a hyperbolic von Dyck group
of signature $\bsdelta=(\delta_1,\delta_2,\delta_3)$
corresponding to an exceptional unimodal singularity.
Let further $V_\Gamma$ be
the minimal smooth normal-crossing compactification 
of the minimal good resolution of the surface
$\Vbar_\Gamma$.
The complement
$
 V_\Gamma \setminus (\Vbar_\Gamma \setminus \bszero)
$
has two connected components.
One is $E_0$,
which is the exceptional locus
of the minimal good resolution
of the exceptional unimodal singularity.
The other is $E_\infty$,
which is a tree of $(-2)$-curves
whose dual graph is the $T_{p,q,r}$-graph
shown in \pref{fg:Tpqr}.

Choose homogeneous generators $x$, $y$ and $z$
of $A(\Gamma)$
of degrees $a$, $b$, and $c$ respectively, and
write $A(\Gamma) = \bC[x,y,z] / (f)$,
where $f \in \bC[x,y,z]$ is a homogeneous polynomial
of degree $h$.
Let $w$ be an indeterminant of degree $h-a-b-c$ and
$F$ be a homogeneous element of $\bC[x,y,z,w]$
of degree $h$
such that $F(x,y,z,0) = f(x,y,z)$.
Assume that $F$ does not have a critical point
outside of the origin, and
let $S=\bC[x,y,z,w]/(F)$ be the quotient ring of $\bC[x,y,z,w]$
by the ideal generated by $F$.
Then the stack
$
 \bY
  = \bProj S
  := [(\Spec S \setminus \bszero) /\bCx]
$
is a smooth Deligne-Mumdord stack,
which has the trivial canonical bundle
since $\deg F = \deg x + \deg y + \deg z + \deg w$.
Let $\Ybar$ be the coarse moduli scheme of $\bY$ and
$Y \to \Ybar$ be the minimal resolution of $Y$.
The scheme $Y$ is a smooth compactification
of the Milnor fiber of $f$.
The condition $F(x,y,z,0) = f(x,y,z)$ impies that
the ring $R = S/(w)$ is isomorphic to $A(\Gamma)$,
so that the divisor $\{ w = 0 \} = \bProj R \subset \bY$
at infinity is isomorphic to the orbifold curve $\bX$.
The Calabi-Yau property of $\bY$ and the adjunction formula implies
that the normal bundle of $\bX$ in $\bY$ is isomorphic
to the cotangent bundle of $\bX$.
Since $\bX$ has three orbifold points
of orders $\delta_1$, $\delta_2$, and $\delta_3$,
the coarse moduli scheme $Y$ has simple singularities
of types $A_{\delta_1-1}$, $A_{\delta_2-1}$, and $A_{\delta_3-1}$
at these orbifold points.
It follows that the minimal resolution $Y$ has a configuration
of $(-2)$-curves whose dual graph is
the $T_\bsdelta$-graph
shown in \pref{fg:Tpqr}.
Here, the central node corresponds to the strict transform
of the coarse moduli scheme $X$ of $\bX$,
and three legs comes from resolutions
of simple singularities of type $A$.
The scheme $Y$ is a K3 surface,
whose N\'{e}ron-Severi lattice is generated by these $(-2)$-curves
for very general $F$.
On the other hand,
the Milnor lattice of $f$ is given by
$\That_\bsgamma \cong T_\bsgamma \bot U$,
where $\bsgamma = (\gamma_1,\gamma_2,\gamma_3)$
is the Gabrielov number of $f$.
Since $E_\infty$ is disjoint from the Milnor fiber of $f$,
the transcendental lattice
\begin{align}
 T(Y)
  = \NS(Y)^\bot
  \subset H^2(Y; \bZ)
\end{align}
of $Y$ clearly contains the Milnor lattice
$\That_\bsgamma$.
One can show
that the embedding of the $\That_\bsgamma$-lattice
into the orthogonal lattice of $T_\bsdelta$
in $H^2(Y; \bZ) \cong E_8 \bot E_8 \bot U \bot U \bot U$
is unique and surjective.
It follows that $T(Y) \cong \That_\bsgamma$
for very general $F$.

Let $\fv$ be the defining polynomial
of the weighted homogeneous exceptional unimodal singularity,
which is related to $f$ by the strange duality
so that its Dolgachev and Gabrielov numbers
$\bsdeltav$ and $\bsgammav$ satisfy
$\bsdeltav=\bsgamma$ and $\bsgammav=\bsdelta$.
Then the transcendental lattice and the algebraic lattice
of the corresponding K3 surface $\Yv$ satisfy
\begin{align}
 \NS(Y) \bot U \cong T(\Yv), \qquad
 T(Y) \cong \NS(\Yv) \bot U.
\end{align}

\subsection{Categorifications of strange duality}
 \label{sc:category}

The Grothendieck group $K(\cT)$
of a triangulated category $\cT$
has a structure of a lattice
with respect to the Euler pairing
\begin{align}
 \la X, Y \ra := \sum_{n \in \bZ} (-1)^n \dim \Ext^n(X, Y).
\end{align}
A \emph{categorification} of a lattice
is a triangulated category
whose Grothendieck group is isometric
to the lattice.

Let $f \in \bC[x,y,z]$ be a homogeneous polynomial in three variables
defining an exceptional unimodal singularity, and
$Y$ be a compactification
of the Milnor fiber of $f$
as in \pref{sc:SD_K3}.
%
The \emph{numerical Grothendieck group} $\cN(Y)$
is the quotient of the Grothendieck group $K(Y)$
by the radical of the Euler form.
Riemann-Roch theorem implies that
$\cN(Y)$ is isomorphic to the lattice
$H^0(Y; \bZ) \oplus \NS(Y) \oplus H^4(Y; \bZ)$
of algebraic cycles
equipped with the Mukai pairing
\begin{align} \label{eq:Mukai_pairing}
 ((a_0, a_2, a_4), (b_0, b_2, b_4))
  = a_2 \cdot b_2 - a_0 b_4 - a_4 b_0.
\end{align}
This allows one to think of $D^b \coh Y$
as a categorification of the algebraic lattice,
in a slightly weak sense
in that one considers the numerical Grothendieck group
instead of the Grothendieck group.


On the dual side,
the \emph{Fukaya category} $\Fuk \Yv$
is an $A_\infty$-category
whose objects are Lagrangian submanifolds of $\Yv$ and
whose spaces of morphisms are Lagrangian intersection Floer complexes
\cite{Fukaya_MHACFH,Fukaya-Oh-Ohta-Ono}.
Since the Euler number of the Lagrangian intersection Floer complex
is the algebraic intersection number of the Lagrangian submanifolds,
the numerical Grothendieck group of the Fukaya category
is a sublattice of the transcendental lattice.
Since the transcendental lattice of $\Yv$ is generated by vanishing cycles,
the Fukaya category $\Fuk \Yv$ is a categorification of the Milnor lattice,
and strange duality is a categorification
of (conjectural) homological mirror symmetry \eqref{eq:hms}.


There is an alternative way to categorify strange duality,
using stable derived categories
\cite{Eisenbud_HACI,Buchweitz_MCM,Orlov_TCS}
and Fukaya-Seidel categories
\cite{Seidel_VC, Seidel_VC2, Seidel_PL}.
Let $R = \bC[x,y,z]/(f)$ be the coordinate ring
of a weighted homogeneous exceptional unimodal singularity, and
$\gr R$ be the category of finitely-generated $\bZ$-graded $R$-modules.
A complex of $\bZ$-graded $R$-module is \emph{perfect}
if it is quasi-isomorphic to a bounded complex of
finitely-generated projective modules.
The \emph{stable derived category}
$\Dbsing(\gr R)$
is defined as the quotient category
of the bounded derived category $D^b(\gr R)$
by the full subcategory
consisting of perfect complexes.
On the dual side,
the \emph{Fukaya-Seidel category}
is the $A_\infty$-category
whose objects are vanishing cycles of $\fv$
and whose spaces of morphisms are
Lagrangian intersection Floer complexes.

\begin{conjecture} \label{cj:hms_f}
There is an equivalence
\begin{align} \label{eq:hms_f}
 \Dbsing(\gr R) \cong D^b \cF(\fv)
\end{align}
of triangulated categories.
\end{conjecture}

The left hand side of \eqref{eq:hms_f}
is a categorification of the algebraic lattice
$\cN(Y) \cong \NS(Y) \bot U \cong \That_\bsdelta$,
in the sense that the Grothendieck group of $\Dbsing(\gr R)$
equipped with the symmetrized Euler pairing
$
 (X, Y) = \la X, Y \ra + \la Y, X \ra
$
is isometric to $\cN(Y)$
\cite{Kobayashi-Mase-Ueda_EUS,MR3223358}.
Similarly,
the right hand side of \eqref{eq:hms_f} is a categorification
of the Milnor lattice of $\fv$.
\pref{cj:hms_f} is known
for a disconnected sum of polynomials of type $A$ or $D$
\cite{Futaki-Ueda_BP,Futaki-Ueda_Dn}.
Fukaya category of unimodal singularities
has been studied by Keating,
who in particular has proved homological mirror symmetry
\cite[Theorem 7.1]{1405.0744}
\begin{align}
 D^b \cF(T_{p,q,r}) \cong D^b \coh \bX_{p,r,q}
\end{align}
between $T_{p,q,r}$-singularity
and orbifold rational curves
conjectured by Takahashi
\cite[Conjecture 7.4]{Takahashi_WPL}.

\subsection{Generalizations of strange duality}


\subsubsection{Cusp singularities}

An isolated singularity in dimension 2 is a \emph{cusp singularity}
if the exceptional locus of the minimal resolution is
either a cycle of non-singular rational curves
or a nodal rational curve.
A cusp singularity is a hypersurface singularity
if and only if it is a $T_{p,q,r}$-singularity
with $1/p+1/q+1/r<1$
\cite{MR0472811}.
Let $C_1, \ldots, C_n$ be the irreducible components
of the exceptional locus.
The off-diagonal part of the intersection matrix
$((C_i, C_j))_{i,j=1}^n$
is given by
$C_i \cdot C_j = 1$ for $|i-j|=1$
and $C_i \cdot C_j = 0$ for $|i-j|>1$,
and only the diagonal part
\begin{align}
 b_i := - C_i \cdot C_i
\end{align}
contains a non-trivial information.
The \emph{cycle number}
of a cusp singularity is
the sequence of integers defined by
\begin{align}
  (d_1, \ldots, d_n)
  =
\begin{cases}
 (b_1, \ldots, b_n) & n \ge 2, \\
 (b_1+2) & n=1.
\end{cases}
\end{align}
Let
\begin{align}
 \omega
  = [\overline{b_1,b_2,\ldots,b_n}]
  = b_1-\cfrac{1}{b_2-\cfrac{1}{b_3-\cfrac{1}{\cdots}}}
\end{align}
be the totally-real quadratic irrational number
defined by an infinite continued fraction.
The \emph{dual cycle number} of the cusp singularity
is the sequence $(e_1,\ldots,e_k)$ of integers
defined by the continued fraction expansion
\begin{align}
 1/\omega = [f_1,\ldots,f_s,\overline{e_1,\ldots,e_k}].
\end{align}
Then a certain class of cusp singularities come in pairs
in such a way that the cycle numbers
and the dual cycle numbers are interchanged.

The duality of cusp singularities can be described
in terms of hyperbolic Inoue surfaces,
just as the strange duality of Arnold can be interpreted
in terms of K3 surfaces.
In the Enriques-Kodaira classification of compact complex surfaces,
a surface of class V\!I\!I is a non-K\"{a}hler elliptic surface
with first Betti number 1.
Inoue surfaces form an important class of surfaces of class V\!I\!I,
which can further be divided into three subclasses;
hyperbolic, half, and parabolic.

A \emph{hyperbolic Inoue surface} is constructed as follows.
Let $K$ be a real quadratic field, and
$M$ be a free $\bZ$-module of rank 2 in $K$.
Define 
\begin{align}
 U_+(M) = \{ x \in K \mid x M = M, x>0, x'>0 \},
\end{align}
where $(-)' \colon K \to K$ is the conjugation.
Let $V$ be a subgroup of $U_+(M)$ of finite index.
The group
\begin{align}
 G(M,V)
  = \lc \sqmat{\varepsilon}{m}{0}{1} \in \GL_2(\bR) \relmid
      \varepsilon \in V, m \in M
     \rc
\end{align}
acts on $\bH \times \bC$ by
\begin{align}
 \sqmat{\varepsilon}{m}{0}{1} \colon
  (z_1, z_2) \mapsto (\varepsilon z_1 + m, \varepsilon' z_2 + m').
\end{align}
The quotient space
$Y^\circ(M,V) = (\bH \times \bC)/G(M,V)$
can be compactified to a compact complex surface $\Ybar(M,V)$
by adding two points $\infty$ and $\infty'$,
which are cusp singularities.
The minimal resolution $Y(M,V)$ of $\Ybar(M,V)$
is a surface of class V\!I\!I,
whose second Betti number is the sum of the length
of the chains of rational curves
obtained from the cusp singularities $\infty$ and $\infty'$
\cite{MR0442297}.
Any minimal surface of class V\!I\!I
with two cycles of rational curves
is isomorphic to a hyperbolic Inoue surface
\cite{MR768987}.
A pair of cusp singularities are dual
if they come from two cusps of a hyperbolic Inoue surface
\cite{MR593635}.
The relation between cusp singularities and mirror symmetry
is investigated in \cite{1106.4977}
using the techniques developed in
\cite{MR2181810,Gross-Siebert_RAGCG}.

\subsubsection{Ebeling-Wall duality}

Strange duality is extended in \cite{MR806842}
to encompass certain classes of bimodal singularities,
isolated complete intersection singularities, and
quadrilateral singularities.
Since a well-written review already exists
\cite{MR1709346},
we do not discuss it here.

\subsubsection{Saito duality}

A sequence $W=(a,b,c;h)$ of positive integers
with $h > \max \{ a,b,c \}$ is a \emph{regular weight system}
if
\begin{align}
 \frac{(1-T^{h-a})(1-T^{h-b})(1-T^{h-c})}
  {(1-T^a)(1-T^b-1)(1-T^c)}
\end{align}
is a polynomial in $T$
\cite{MR894306}.
The integers $(a,b,c)$ and $h$ are called
\emph{weights} and the \emph{Coxeter number}
respectively.
For a regular weight system $W$,
define a sequence $(m_1, \ldots, m_\mu)$
of not necessarily distinct integers by
\begin{align}
 \chi_W(T) &:= T^{\epsilon} \,
  \frac{(1-T^{h-a})(1-T^{h-b})(1-T^{h-c})}
   {(1-T^a)(1-T^b-1)(1-T^c)}
 = T^{m_1} + \cdots + T^{m_\mu},
\end{align}
where
\begin{align}
 \epsilon
  = \epsilon_W
  := a+b+c-h
\end{align}
is the \emph{minimal exponent} of $W$, and
\begin{align}
 \mu
  = \mu_W
  := \chi_W(1)
  =(h-a)(h-b)(h-c)/abc
\end{align}
is the \emph{rank} of $W$.
The \emph{characteristic polynomial} of $W$ is defined by
\begin{align}
 \varphi_W(\lambda)
  = \prod_{i=1}^\mu (\lambda - \omega_i),
\end{align}
where
\begin{align}
 \omega_i = \exp(2 \pi \sqrt{-1} m_i/h),
  \qquad i=1, \ldots, \mu.
\end{align}
One can write
\begin{align}
 \varphi(\lambda)
  = \prod_{i|h} (\lambda^i-1)^{e(i)}
\end{align}
for some sequence $(e(i))_{i|h}$
of integers.
A pair $(W, W^*)$ of regular weight systems
with the identical Coxeter number
are said to be \emph{dual}
in the sense of Saito
\cite{MR1734136}
if their characteristic polynomials
\begin{align}
 \varphi_W(\lambda)
  = \prod_{i|h} (\lambda^i-1)^{e(i)},
 \qquad
 \varphi_{W^*}(\lambda)
  = \prod_{i|h} (\lambda^i-1)^{e^*(i)},
\end{align}
satisfy
\begin{align}
 e(i) + e^*(h/i) = 0
\end{align}
for all $i|h$.

A weight system $W=(a,b,c;h)$ is regular
if and only if there exists a weighted homogeneous polynomial
$f \in \bC[x,y,z]$
of weight $\deg(x,y,z;f) = (a,b,c;h)$
with isolated critical point at the origin.
The rank of a regular weight system is the Milnor number
$\dim \bC[x,y,z]/(\partial_x f, \partial_y f, \partial_z f)$
of the singularity,
and the minimal exponent is the Gorenstein parameter
of the coordinate ring $\bC[x,y,z]/(f)$.
The characteristic polynomial
is the characteristic polynomial
for the Milnor monodromy
acting on the middle-dimensional homology group
of the Milnor fiber.

The duality of regular weight systems
generalizes the strange duality of exceptional unimodal singularities.
See
\cite{MR1653033} for a summary
of the theory of regular weight systems, and
\cite{MR1734136} for a thorough treatment.
The relation
between the duality of regular weight systems
and mirror symmetry
is discussed in
\cite{Kobayashi_DW,
Takahashi_DRWS,
MR2278769}.

\subsubsection{Kobayashi duality}

A \emph{weight system}
in the sense of Kobayashi
\cite{Kobayashi_DW}
is a sequence
$W=(a_1,\ldots,a_n;h)$ of positive integers
satisfying $h \in \bN a_1 + \cdots + \bN a_n$.
A weight system is \emph{reduced}
if $\gcd(a_1,\ldots,a_n;h)=1$.
A \emph{weighted magic square}
for a pair $(W, W^*)$ of weight systems
is an $n \times n$ matrix $C$
with non-negative integer entries
satisfying
\begin{align}
 C
\begin{pmatrix}
 a_1 \\ \vdots \\ a_n
\end{pmatrix}
 =
\begin{pmatrix}
 h \\ \vdots \\ h
\end{pmatrix}
\quad
\text{and}
\quad
\begin{pmatrix}
 a_1^* & \cdots & a_n^*
\end{pmatrix}
 C
 =
\begin{pmatrix}
 h^* & \cdots & h^*
\end{pmatrix}.
\end{align}
A weighted magic square $C$ is \emph{primitive}
if $|\det C| = h = h^*$.
The weight systems $W$ and $W^*$ are \emph{dual}
in the sense of Kobayashi
\cite{Kobayashi_DW}
if there exists a primitive weighted magic square
for $(W, W^*)$.
This also gives a generalization
of the strange duality of exceptional unimodal singularities.
This notion of duality is further investigated in
\cite{MR2278769}.

\section{Berglund-H\"{u}bsch transpose}
 \label{sc:BH}

A polynomial $f \in \bC[x_1, \ldots, x_n]$ is {\em invertible}
if the following two conditions are satisfied:
\begin{enumerate}
 \item
There is an integer matrix
$
 A = (a_{ij})_{i, j=1}^n
$
such that
\begin{align}
 f = \sum_{i=1}^n \prod_{j=1}^n x_j^{a_{ij}}.
\end{align}
 \item
$f$ has an isolated critical point at the origin.
\end{enumerate}
The second condition implies that
the matrix $A$ has a non-zero determinant.
It follows that any invertible polynomial is weighted homogeneous,
and the corresponding reduced weight system
\begin{align}
 (a_1, \ldots, a_n; h) := \deg(x_1, \ldots, x_n; f)
\end{align}
is determined uniquely.

An invertible polynomial is a Sebastiani-Thom
(or decoupled) sum of polynomials
of the following three types
\cite{MR1188500}:
\begin{itemize}
 \item
Fermat:$x^p.$
 \item
Chain:
$
 x_1^{p_1} x_2 + x_2^{p_2} x_3
  + \cdots + x_{n-1}^{p_{n-1}} x_n + x_n^{p_n}.
$
 \item
Loop:
$
 x_1^{p_1} x_2 + x_2^{p_2} x_3
  + \cdots + x_{n-1}^{p_{n-1}} x_n + x_n^{p_n} x_1.
$
\end{itemize}

A \emph{Landau-Ginzburg model} is a pair
$(V, f)$
of an algebraic variety $V$
and a holomorphic function $f \colon V \to \bC$ on $V$.
An invertible polynomial
$
 f
$
gives an example $(\bC^n, f)$
of a Landau-Ginzburg model.
 
An invertible polynomial is naturally graded
by the abelian group $L$ of rank one
generated by $n+1$ elements
$\vecx_1, \ldots, \vecx_n$, and $\vecc$
with $n$ relations
\begin{align}
 a_{i1} \vecx_1 + \cdots + a_{in} \vecx_n = \vecc,
  \qquad i=1, \ldots, n.
\end{align}
The group $L$ is the group of characters of the group $K$
defined by
\begin{align}
 K = \{ (\alpha_1, \ldots, \alpha_n) \in (\bCx)^n \mid
  \alpha_1^{a_{11}} \cdots \alpha_n^{a_{1n}}
   = \cdots
   = \alpha_1^{a_{n1}} \cdots \alpha_n^{a_{nn}} \}.
\end{align}
The \emph{group of maximal diagonal symmetries}
is the kernel $\Gmax$ of the map
\begin{align}
\begin{array}{ccc}
 K & \to & \bCx \\
 \vin & & \vin \\
 (\alpha_1, \ldots, \alpha_n) & \mapsto &
  \alpha_1 \cdots \alpha_n,
\end{array}
\end{align}
so that there is an exact sequence
\begin{align}
 1 \to \Gmax \to K \to \bCx \to 1.
\end{align}
This exact sequence induces an exact sequence
\begin{align}
 1 \to \bZ \to L \to \Gmax^\vee \to 1
\end{align}
of the corresponding groups of characters,
where
\begin{align}
 \Gmax^\vee = \Hom(\Gmax, \bCx)
\end{align}
is (non-canonically) isomorphic to $\Gmax$.
Write
\begin{align}
 A^{-1} =
\begin{pmatrix}
 \varphi_1^{(1)} & \varphi_1^{(2)} & \cdots & \varphi_1^{(n)} \\
 \varphi_2^{(1)} & \varphi_2^{(2)} & \cdots & \varphi_2^{(n)} \\
 \vdots & \vdots & \ddots & \vdots \\
 \varphi_n^{(1)} & \varphi_n^{(2)} & \cdots & \varphi_n^{(n)}
\end{pmatrix}.
\end{align}
Then the group $\Gmax$ is generated by
\begin{align}
 \rho_k = \lb \exp \lb 2 \pi \sqrt{-1} \varphi_1^{(k)} \rb, \ldots,
  \exp \lb 2 \pi \sqrt{-1} \varphi_n^{(k)} \rb \rb,
 \qquad k=1, \ldots, n.
\end{align}

A \emph{Landau-Ginzburg orbifold} is a pair
$((V, f), G)$ of a Landau-Ginzburg model $(V, f)$
and a finite group $G$ acting on $V$
in such a way that the function $f \colon V \to \bC$
is $G$-invariant.
An invertible polynomial $f$ and a subgroup $G$ of $\Gmax$
gives an example
$((f, \bC^n), G)$
of a Landau-Ginzburg orbifold.

The \emph{transpose} of $f$ is the invertible polynomial
\begin{align}
 \fv =  \sum_{i=1}^n \prod_{j=1}^n x_j^{a_{ji}}.
\end{align}
Note that the exponent matrix $\Av$
is the transpose matrix of $A$.
The group $\Gvmax$ of maximal diagonal symmetries of $\fv$
is generated by the column vectors
$\rhov_i$ of $\Av^{-1}$.
The \emph{transpose} of a subgroup $G \subset \Gmax$
is defined in \cite{MR1310310,Krawitz} as
\begin{align}
 \Gv = \lc \prod_{i=1}^n \rhov_i^{r_i} \relmid
  \begin{pmatrix} r_1 & \cdots & r_n \end{pmatrix}
  A^{-1}
  \begin{pmatrix} a_1 \\ \vdots \\ a_n \end{pmatrix}
  \in \bZ \text{ for all }
  \prod_{i=1}^n \rho_i^{a_i} \in G \rc.
\end{align}
The \emph{transpose} of a pair $(f, G)$
of an invertible polynomial $f$
and a subgroup $G$ of the group $\Gmax$
of maximal diagonal symmetries
is defined as $(\fv, \Gv)$.
In particular,
the transpose of the pair
$(f, 1)$
of an invertible polynomial
and the trivial group
is the pair
$(\fv, \Gvmax)$
of the transpose polynomial
and the group of maximal diagonal symmetries.

Transposition is introduced in \cite{Berglund-Hubsch}
as a generalization of the orbifold mirror construction
\cite{Greene-Plesser}.
Interpretation of strange duality
as mirror symmetry for Landau-Ginzburg orbifolds
goes back at least to
\cite{MR1079045,MR1104265}.
There has been a renewed interest
in transposition mirror construction recently,
coming in part from the development
of cohomological field theories
associated with invertible polynomials
\cite{MR3043578, 1105.2903}.
In particular,
classical mirror symmetry for exceptional unimodal singularities
is proved in \cite{1405.4530}.
Transposition mirror construction
as a generalization of strange duality
is explored in \cite{Ebeling-Takahashi_SDWHP}.

To relate mirror symmetry of Landau-Ginzburg orbifolds
to mirror symmetry of Calabi-Yau manifolds,
one needs the correspondence
between Calabi-Yau manifolds and Landau-Ginzburg orbifolds.
Let $f$ be an invertible polynomial and
$G$ be a subgroup of $\Gmax$.
Put
\begin{align}
 \varphi_i = \varphi_i^{(1)} + \cdots + \varphi_i^{(n)},
  \qquad i = 1, \ldots, n,
\end{align}
and define a homomorphism
$
 \varphi \colon \bCx \to K
$
by
\begin{align}
 \varphi(\alpha)
  = \lb \alpha^{h \varphi_1}, \ldots, \alpha^{h \varphi_n} \rb,
\end{align}
where $h$ is the smallest positive integer
such that $a_i := h \varphi_i \in \bZ$
for $i=1,\ldots,n$.
Then $\varphi$ is injective and
one has an exact sequence
\begin{align}
 1 \to \bCx \xto{\varphi} K \to \Gbarmax \to 1,
\end{align}
where $\Gbarmax := \coker \varphi$.
The polynomial $f$ is weighted homogeneous
with respect to the weight system $(a_1, \ldots, a_n; h)$.
Since $f$ has an isolated singularity at the origin,
the weighted projective hypersurface
\begin{align}
 \bY = \{ [x_1:\cdots:x_n] \in \bP(a_1, \ldots, a_n)
  \mid f(x_1,\ldots,x_n)=0 \}
\end{align}
is a smooth Deligne-Mumford stack.
It is Calabi-Yau if
$
 a_1+\cdots+a_n=h.
$
The intersection $\Image \varphi \cap \Gmax$ is generated by
\begin{align}
 J = \lb \exp \lb 2 \pi \sqrt{-1} \varphi_1 \rb, \ldots,
  \exp \lb 2 \pi \sqrt{-1} \varphi_n \rb \rb.
\end{align}
Assume $J \in G$ and let $\Gbar = G / \la J \ra$
be the image of $G$ in $\Gbarmax$.
Then $\Gbar$ acts naturally on $\bY$,
and one can form the quotient stack
$[\bY/\Gbar]$.
The \emph{Calabi-Yau/Landau-Ginzburg correspondence},
or the CY/LG correspondence for short,
is an idea
which goes back at least to \cite{MR929596},
that the Landau-Ginzburg orbifold $(f, G)$
is `dual' to the orbifold $[\bY/\Gbar]$.
On the symplectic side,
this is proved
at the level of topological mirror symmetry
in \cite{MR2672282}.
On the complex side,
the equivalence of derived categories is established
in \cite{Orlov_DCCSTCS}.
The relation between categorical equivalence
and analytic continuation of periods is discussed
in \cite{1201.0813}.

Transposition mirror construction
gives a candidate for the mirror of a Landau-Ginzburg orbifold,
which produces a candidate for the mirror of a Calabi-Yau manifold
by the CY/LG correspondence.
For K3 surfaces,
this construction is known to give mirror pairs
in the sense of \pref{df:dolgachev} below
for an invertible polynomial of the form
$
 x^p + f(y,z,w)
$
for $p=2$ in \cite{Artebani-Boissiere-Sarti}
and for prime $p$ in \cite{1211.2172}.


\section{Aspinwall-Morrison mirrors}
 \label{sc:am}

Let $Y$ be a K3 surface.
The {\em extended K3 lattice} is the free abelian group
$
 H^*(Y; \bZ)
$
equipped with the Mukai pairing
\eqref{eq:Mukai_pairing}.
The complex structure of $Y$ is determined
by the class
$
 [\Omega] \in \bP(H^2(Y, \bC))
$
of the holomorphic 2-form
$
 \Omega \in H^{2,0}(Y)
$
satisfying the Hodge-Riemann bilinear relations
\begin{align}
 (\Omega, \Omega) = 0, \qquad
 (\Omega, \Omegabar) > 0.
\end{align}
A {\em complexified K\"{a}hler structure} on $Y$
is an element $\mho \in H^*(Y, \bC)$ of the form
\begin{align*}
 \mho
 &= \exp(B + \sqrt{-1} \omega) \\
 &= \lb 1, B + \sqrt{-1} \omega,
      \frac{1}{2} (B + \sqrt{-1} \omega)^2 \rb \\
 &\in H^0(Y;\bC) \oplus H^2(Y;\bC) \oplus H^4(Y;\bC)
\end{align*}
satisfying
\begin{align}
 (\mho, \mho) = 0, \qquad
 (\mho, \mhobar) > 0
\end{align}
and
\begin{align}
 (\mho, \Omega) = 0, \qquad
 (\mho, \Omegabar) = 0.
\end{align}
The class $\omega \in H^2(Y, \bR)$ is
(a slight generalization of) the K\"{a}hler class,
which satisfies $\omega \in H^2(Y, \bR) \cap H^{1,1}(Y)$ and
$\omega \cdot \omega > 0$.
The class $B \in H^2(Y, \bR)$ is called
the {\em $B$-field},
which satisfies $B \in H^2(Y, \bR) \cap H^{1,1}(Y)$.
Note that $\omega$ and $B$ are determined
by the class
$
 [\mho] \in \bP(H^*(Y, \bC)).
$

Let $\Yv$ be another K3 surface and
$(\Omegav, \mhov)$ be a pair of 2-forms
satisfying the above conditions.
The following definition is due to
Aspinwall and Morrison
\cite{Aspinwall-Morrison_STKS}:

\begin{definition}
A pair $((Y,(\Omega, \mho)), (\Yv,(\Omegav, \mhov)))$
of K3 surfaces equipped with complexified K\"{a}hler structures
is a {\em mirror pair} if
there is an isometry of extended K3 lattices
\begin{align}
 \varphi: H^*(Y; \bZ) \simto H^*(\Yv; \bZ)
\end{align}
such that
\begin{align}
 (\varphi(\Omega), \varphi(\mho)) = (\mhov, \Omegav).
\end{align}
\end{definition}

\section{Dolgachev mirrors}
 \label{sc:Dolgachev}

Let $Y$ be a K3 surface.
The {\em K3 lattice} is the free abelian group
$H^2(Y, \bZ)$
equipped with the intersection form.
It has rank 22 and signature (3, 19),
and isomorphic to
\begin{align}
 L = E_8 \bot E_8 \bot U \bot U \bot U
\end{align}
as an abstract lattice.
Here $E_8$ is the negative-definite even unimodular lattice
of type $E_8$ and
$U$ is the even unimodular indefinite lattice
of rank two.
For a K3 surface $Y$,
set
\begin{align}
 \Delta(Y) = \{ \delta \in \Pic(Y) \mid (\delta, \delta) = -2 \}.
\end{align}
Let $\cL$ be a line bundle
such that $[\cL] = \delta \in \Pic(Y)$.
Riemann-Roch theorem gives
\begin{align}
 h^0(\scL) + h^0(\scL^\vee)
  \ge 2 + \frac{1}{2}(\delta, \delta),
\end{align}
so that either $\scL$ or $\scL^\vee$ has a non-trivial section,
and hence either $\delta$ or $- \delta$ is effective;
\begin{align}
 \Delta(Y) &= \Delta(Y)^+ \amalg \Delta(Y)^-, \\
 \Delta(Y)^+ &= \{ \delta \in \Delta(Y)
  \mid \delta \text{ is effective} \}, \\
 \Delta(Y)^- &= - \Delta(Y)^+.
\end{align}
The subgroup $W(Y) \subset O(L)$ generated by reflections
with respect to elements in $\Delta(Y)$ acts
properly discontinuously on the connected component
$V^+$ of
\begin{align}
 V(Y) = \{ x \in H^{1,1}(Y) \cap H^2(Y, \bR)
  \mid (x, x) > 0 \}
\end{align}
containing the K\"{a}hler class.
The fundamental domain is given by
\begin{align}
 C(Y) = \{ x \in V(Y)^+ \mid (x, \delta) \ge 0 \text{ for any }
  \delta \in \Delta(Y)^+ \},
\end{align}
and the K\"{a}hler cone is given
(cf.~e.g.~\cite[Corollary VIII.3.9]{Barth-Hulek-Peters-Van_de_Ven})
by
\begin{align}
 C(Y)^+ = \{ x \in V(Y)^+ \mid (x, \delta) > 0 \text{ for any }
  \delta \in \Delta(Y)^+ \}.
\end{align}
Recall that
\begin{align}
 \Pic(Y) = H^{1,1}(Y) \cap H^2(Y; \bZ)
\end{align}
by the Lefschetz theorem.
Set
\begin{align}
 \Pic(Y)^+ = C(Y) \cap H^2(Y; \bZ), \\
 \Pic(Y)^{++} = C(Y)^+ \cap H^2(Y; \bZ).
\end{align}

Let $M$ be an even non-degenerate lattice
of signature $(1, t)$
where $0 \le t \le 19$.
Choose one of two connected components of
\begin{align}
 V(M) = \lc x \in M_\bR \relmid (x,x) > 0 \rc
\end{align}
and call it $V(M)^+$.
Choose a subset $\Delta(M)^+$ of
\begin{align}
 \Delta(M) = \lc \delta \in M \relmid (\delta, \delta) = -2 \rc
\end{align}
such that
\begin{enumerate}
 \item
$\Delta(M) = \Delta(M)^+ \amalg \Delta(M)^-$
where $\Delta(M)^- = \lc - \delta \relmid \delta \in \Delta(M)^+ \rc$, and
 \item
$\Delta(M)^+$ is closed under addition (but not subtraction).
\end{enumerate}
Define
\begin{align}
 C(M)^+ = \lc h \in V(M)^+ \cap M \relmid
  (h, \delta) > 0 \text{ for all } \delta \in \Delta(M)^+ \rc.
\end{align}

The notion of a lattice-polarized K3 surface
is introduced by Nikulin \cite{MR544937}.
The following definition is taken from
\cite[Section 1]{Dolgachev_MSK3}:

\begin{definition}
An {\em $M$-polarized K3 surface} is a pair $(Y, j)$
where $Y$ is a K3 surface and
$
 j : M \hookrightarrow \Pic(Y)
$
is a primitive lattice embedding.
An {\em isomorphism} of $M$-polarized K3 surfaces
$(Y, j)$ and $(Y', j')$ is an isomorphism
$f : Y \to Y'$ of K3 surfaces
such that $j = f^* \circ j'$.
An $M$-polarizaed K3 surface is {\em pseudo-ample} if
\begin{align}
 j(C(M)^+) \cap \Pic(Y)^+ \ne \emptyset,
\end{align}
and {\em ample} if
\begin{align}
 j(C(M)^+) \cap \Pic(Y)^{++} \ne \emptyset.
\end{align}
\end{definition}

Assume that for any two primitive embeddings
$
 \iota_1, \iota_2 \colon M \hookrightarrow L,
$
there exists an isometry
$
 \sigma \colon L \to L
$
such that
$
 \sigma \circ \iota_1 = \iota_2.
$
Fix a primitive lattice embedding
$
 i_M : M \hookrightarrow L
$
and let
$
 N = M^\bot
$
be the orthogonal complement.
The period domain
\begin{align}
 \cD(M)
  = \{ [\Omega] \in \bP(N \otimes \bC)
  \mid (\Omega, \Omega) = 0, \ (\Omega, \Omegabar) > 0 \}
\end{align}
can be identified with the symmetric homogeneous space
$
 O(2, 19-t) / SO(2) \times O(19-t)
$
of oriented positive-definite 2-planes in $N_\bR$.
It consists of two onnected components
$\cD^+(M)$ and $\cD^-(M)$,
each of which is isomorphic to a bounded Hermitian domain
of type IV. Set
\begin{align}
 \Gamma(M) = \{ \sigma \in O(L) \mid \sigma(m) = m
  \text{ for any } m \in M \}
\end{align}
and $\Gamma_M$ be its image
under the natural injective homomorphism
$
 \Gamma(M) \hookrightarrow O(N).
$
Global Torelli theorem and surjectivity of the period map
show that
$
 \cD_0^+(M) / \Gamma_M
$
is the coarse moduli space of ample $M$-polarized K3 surfaces,
where
\begin{align}
 \cD_0^+(M) = \cD^+(M) \setminus \lb \bigcup_{\delta \in \Delta(N)}
  \delta^\bot \rb
\end{align}
is the complement of reflection hyperplanes
\begin{align}
 \delta^\bot = \{ z \in \cD^+(M) \mid (z, \delta) = 0 \}.
\end{align}
The closure of $\cD^+(M)$
in the {\em compact dual}
\begin{align}
 \cDv(M) = \{ [\Omega] \in \bP(N \otimes \bC) \mid
  (\Omega, \Omega) = 0 \}
\end{align}
of the period domain
is denoted by
$\cDbar(M)$.
Its topological boundary is given by
$$
 \cDbar(M) \setminus \cD^+(M)
  = \bigcup_{I \text{ : isotropic subspace of $M_\bR$}}
      \bP(I_{\bC}) \cap \cDbar(M).
$$
Since the signature of $M$ is $(2, 19-t)$,
one either has $\rank I = 1$ or 2,
so that $\bP(I_{\bC}) \cap \cDbar(M)$ is
one point or isomorphic to the upper half plane.
The boundary component is {\em rational}
if $I$ is defined over $\bQ$.
The {\em Satake-Baily-Borel compactification}
is defined by
\begin{align}
 \lb \cD^+(M) \cup
       \bigcup_{I \text{ : rational}}
         \bP(I_{\bC}) \cap \cDbar(M)
      \rb / \Gamma_M.
\end{align}

For a vector $f$ in a lattice $S$,
the positive integer $\sdiv f$ is defined as
the greatest common divisor of
$(f, g) \in \bZ$
for all $g \in S$.
An isotropic vector $f$ is called {\em $m$-admissible}
if $\sdiv f = m$ and there exists another isotropic vector $g$
such that $(f, g) = m$ and $\sdiv g = m$.
If $M^\bot$ has an $m$-admissible vector $f$,
then one has $M^\bot = U(m) \bot \Mv$,
where $U(m)$ is the lattice generated by $f$ and $g$.

\begin{definition}[{Dolgachev \cite[Section 6]{Dolgachev_MSK3}}]
 \label{df:dolgachev}
The {\em mirror moduli space} is the moduli space
$
 D_{\Mv}^\circ
$
of ample $\Mv$-polarized K3 surfaces.
\end{definition}

The case $m = 1$ is of particular interest.
Two sublattices $M$ and $\Mv$
of the K3 lattice $L$ are said to be
{\em K3-dual} if one has
$$
 M^\bot = \Mv \bot U
$$
for a unimodular hyperbolic plane $U$ in $L$.
Nikulin \cite{Nikulin_ISBF} has shown that
two hyperbolic lattices $M$ and $\Mv$ are K3-dual
if $\rank M + \rank \Mv = 20$ and
the discriminant group
$A(M) = M^* / M$ equipped with
the discriminant form $q_M$
is isomorphic to $A(\Mv)$
equipped with $- q_{\Mv}$.
This duality has been investigated by Belcastro
\cite{Belcastro_PLFK3S}
in the case of K3 surfaces associated with weighted projective spaces.
There are 95 weights
where the minimal model
of a general anti-canonical hypersurface
of the corresponding weighted projective space
gives a K3 surface
\cite{Yonemura_HSKS}.
Belcastro computed the Picard lattice
of very general K3 surfaces obtained
in this way,
and discovered that some of them are K3-dual
with each other.
The relation between Belcastro duality
and Kobayashi duality is studied in \cite{MR2278769}.

\section{Batyrev mirrors}
 \label{sc:Batyrev}

Let $Y$ be a smooth anti-canonical hypersurface
in a smooth toric weak Fano 3-fold $X$.
Here, a {\em weak Fano} manifold is a projective manifold
whose anti-canonical bundle is nef and big.
It follows from the adjunction formula
that $Y$ has the trivial canonical bundle.
Lefschetz hyperplane theorem shows that $Y$ is simply-connected,
so that $Y$ is a K3 surface.
Batyrev \cite{Batyrev_DPMS}
introduced the following construction
of a candidate for the mirror of $Y$.

Let $\bT \subset X$ be the dense torus
and $\bsM = \Hom(\bT, \bCx)$ be the group of characters
so that $\bT = \Spec \bC[\bsM]$.
Let further
\begin{align}
 \Delta = \Conv \{ m \in \bsM \mid x^m \in H^0(\scO_X(-K_X)) \}
\end{align}
be the Newton polytope
of the anti-canonical bundle of $X$,
where the canonical divisor $K_X$ is
the sum of all prime toric divisors, and
$x^m \in \bC[\bsM]$ is
the rational function on $X$
corresponding to $m \in \bsM$.
Since $X$ is weak Fano,
the polytope $\Delta$ is reflexive,
i.e., the polar dual polytope
\begin{align}
 \Deltav = \{ n \in \bsN \mid \la n, m \ra \ge -1
  \text{ for any } m \in \Delta \}
\end{align}
is also a lattice polytope
(i.e., the convex hull of a finite subset of $\bsN$).
Here $\bsN = \Hom(\bsM, \bZ)$
is the group of one-parameter subgroups of $\bT$.

Recall that the {\em fan polytope} of a toric variety
is the convex hull of primitive generators
of one-dimensional cones of the fan.
Let $\Sigmav$ be any unimodular fan
in $\bsM_\bR = \bsM \otimes \bR$
whose fan polytope coincides with $\Delta$.
Note that the fan $\Sigmav$ lives in $\bsM_\bR$,
whereas one usually considers a fan inside $\bsN_\bR$.
This comes from the fact that we are working on the mirror side,
where the roles of the torus and the dual torus are interchanged.
Let $\Xv$ be the toric variety
associated with the fan $\Sigmav$.
The toric variety $\Xv$ is weak Fano
since $\Delta$ is reflexive,
and a general member
$\Yv \in \left| - K_{\Xv} \right|$
of the anti-canonical linear system
is a smooth K3 surface.

Let $\Pictor(Y)$ be the sublattice
of the Picard lattice $\Pic(Y)$
generated by restrictions of the toric divisors of $X$.
Some of the restrictions of the toric divisors may be reducible,
and we write the sublattice of $\Pic(Y)$
generated by irreducible components of restrictions of toric divisors
as $\Piccor(Y)$
following \cite{Whitcher_RP}.
The classification by
Kreuzer and Skarke
\cite{MR1663339}
shows that 
there are 4319 reflexive polytopes in dimension 3.
Rohsiepe \cite{Rohsiepe_LPTK3S}
computed $\Pictor(Y)$ and $\Piccor(Y)$
for all of them,
and showed that $\Piccor(Y)$ and $\Pictor(\Yv)$
are K3-dual.

\section{Classical mirror symmetry}
 \label{sc:classical}

Let $\bsN \cong \bZ^3$ be a free abelian group of rank three and
$\bsM = \Hom(\bsN, \bZ)$ be the dual group.
Let further $(\Sigma, \Sigmav)$ be a pair
of unimodular fans
in $\bsN_\bR = \bsN \otimes \bR$ and
$\bsM_\bR = \bsM \otimes \bR$
whose fan polytopes $(\Delta, \Deltav)$
are polar dual to each other.
The set of generators
of one-dimensional cones of the fan $\Sigma$
will be denoted by $\{ b_1, \ldots, b_m \} \subset \bsN$.
One has the {\em fan sequence}
\begin{align}
 0 \to \bL \to \bZ^m \xto{\beta} \bsN \to 0
\end{align}
and the {\em divisor sequence}
\begin{align}
 0 \to \bsM \xto{\beta^*} (\bZ^m)^* \xto{} \bL^* \to 0,
\end{align}
where $\beta$ sends the $i$-th coordinate vector to $b_i$ and
\begin{align}
 \Pic(X) \cong H^2(X; \bZ) \cong \bL^*.
\end{align}
Set $\scM = \bL^* \otimes \bCx$
and $\bTv = \bsM \otimes \bCx$
so that one has the exact sequence
\begin{align}
 1 \to \bTv \to (\bCx)^m \to \scM \to 1.
\end{align}
The uncompactified mirror $\Yv^\circ_\alpha$
of a smooth anticanonical hypersurface $Y \subset X$
is defined by
\begin{align}
 \Yv^\circ_\alpha = \{ y \in \bTv \mid
  W_\alpha(y) = \sum_{i=1}^{m} \alpha_i y^{b_i} = 1 \}
\end{align}
where $\alpha = (\alpha_1, \ldots, \alpha_m) \in (\bCx)^m$.
We will study symplectic geometry of $Y$,
so that the defining equation of $Y$ is irrelevant.
The closure $\Yv_\alpha$ of $\Yv^\circ_\alpha$ in $\Xv$ for general $\alpha$
is a smooth anti-canonical K3 hypersurface,
which is the compact mirror of $Y$.
We say that $\Yv_\alpha$ is $\Sigmav$-regular
if the intersection of $\Yv_\alpha$
with any toys orbit of $\Xv$ is a smooth subvariety of codimension one.
Let $(\bCx)^m_\reg$ be the set of $\alpha \in (\bCx)^m$
such that $\Yv_\alpha$ is $\Sigmav$-regular, and
$
 \widetilde{\varphiv} : \widetilde{\fYv} \to (\bCx)^m_\reg
$
be the second projection from
\begin{align}
 \widetilde{\fYv} = \{ (y, \alpha) \in \Xv \times (\bCx)^m_\reg
   \mid W_\alpha(y) = 1 \}.
\end{align}
The quotient of the family
$
 \widetilde{\varphiv} : \widetilde{\fYv} \to (\bCx)^m_\reg
$
by the free $\bTv$-action
\begin{align}
 t \cdot (y, (\alpha_1, \ldots, \alpha_m))
  = (t^{-1} y, (t^{b_1} \alpha_1, \ldots, t^{b_m} \alpha_m))
\end{align}
will be denoted by
$
 \varphiv : \fYv \to \scM_\reg,
$
where $\scM_\reg = (\bCx)^m_\reg / \bTv$.
The \emph{residue part} of $H^2(\Yv_\alpha;\bC)$
is defined as the image of the residue map:
\begin{align}
 H^2_\res(\Yv_\alpha;\bC)
  := \Image(\Res \colon H^0(\Xv, \Omega^3_\Xv(*\Yv_\alpha))
   \to H^2(\Yv_\alpha;\bC)).
\end{align}
One can show
\cite[Section 6.3]{Iritani_QCP} that
$H^2_\res(\Yv_\alpha;\bC)$ can be identified
with the lowest weight component
$W_2(H^2(\Yv_\alpha^\circ;\bC))$
of the mixed Hodge structure on $H^2(\Yv_\alpha^\circ;\bC)$,
and hence comes naturally with a $\bQ$-Hodge structure
of weight 2.
The {\em residual B-model VHS}
$
 (\scrH_B, \nabla^B, \HBQ, \scrF_B^\bullet, Q_B)
$
on $\Uv$ consists of
\begin{itemize}
 \item
the locally-free subsheaf $\scrH_B$
of $(R^2 \varphiv_* \bC_\fYv) \otimes \cO_{\cM_\reg}$
whose fiber at $[\alpha]$ is
$H^2_\res(\Yv_\alpha;\bC)$,
 \item
the Gauss-Manin connection $\nabla^B$ on $\scrH_B$,
 \item
the rational structure $\HBQ \subset \Ker \nabla^B$ explained above,
 \item
the standard Hodge filtration
$
 \scrF_{B,[\alpha]}^p = \bigoplus_{j \ge p} H_\res^{j,2-j}(\Yv_\alpha;\bC),
$
and
 \item
the intersection form
\begin{align}
 Q_B(\omega_1, \omega_2)
  = \int_{\Yv_\alpha} \omega_1 \cup \omega_2.
\end{align}
\end{itemize}
The composition
\begin{align}
 H_3(\bTv, \Yv_\alpha^\circ; \bC)
  \xto{\ \partial \ } H_2(\Yv_\alpha^\circ; \bC)
  \to H_2(\Yv_\alpha; \bC)
  \xto{\ \PD \ } H^2(\Yv_\alpha; \bC)
\end{align}
gives a surjection
$
 \VC \colon H_3(\bTv, \Yv_\alpha^\circ; \bC) \to H^2_\res(\Yv_\alpha; \bC)
$
called the \emph{vanishing cycle map}
.
The image of $H_n(\bTv, \Yv_\alpha^\circ; \bZ)$
by the vanishing cycle map
defines the \emph{vanishing cycle integral structure}
$\HBZ^\vc \subset \HBQ$
on the residual B-model VHS
\cite[Definition 6.7]{Iritani_QCP}.

On the A-model side, let
\begin{align}
 H^\bullet_\amb(Y; \bC)
  = \Image(\iota^* : H^\bullet(X; \bC) \to H^\bullet(Y; \bC))
\end{align}
be the subspace of $H^\bullet(Y; \bC)$
coming from the cohomology classes of the ambient toric variety,
and set
\begin{align}
 U = \{ \tau \in H^2_\amb(Y; \bC) \mid
   \Re \la \tau, d \ra \le -M \text{ for any non-zero }
   d \in \Eff(Y) \}
\end{align}
for some sufficiently large $M$.
Here $\Eff(Y)$ is the semigroup of effective curves.
This open subset $U$ is considered
as a neighborhood of the large radius limit point.
Choose an integral basis $p_1, \ldots, p_r$ of $\Pic X$
such that each $p_i$ is nef, and
let $(\tau^i)_{i=1}^r$ be the dual coordinate on $H^2_\amb(Y; \bC)$;
$\tau = \sum_{i=1}^r \tau^i p_i$.

The {\em ambient A-model VHS}
$(\HA, \nabla^A, \scrF_A^\bullet, Q_A)$
\cite[Definition 6.2]{Iritani_QCP}
consists of
\begin{itemize}
 \item
the locally free sheaf
$
 \HA = H_\amb^\bullet(Y) \otimes \scO_U,
$
 \item
the Dubrovin connection
$
 \nabla^A = d + \sum_{i=1}^r (p_i \circ_\tau) \, d \tau^i
  \colon \HA \to \HA \otimes \Omega_{U}^1,
$
 \item
the Hodge filtration
$
 \scrF_A^p = H_\amb^{\le 4-2p}(Y) \otimes \scO_U,
$
and
 \item
the symmetric pairing
$
 Q_A : \HA \otimes \HA \to \scO_U,
$
\ 
$
 (\alpha, \beta) \mapsto (2 \pi \sqrt{-1})^2
  \lb (-1)^{\frac{\deg}{2}} \alpha, \beta \rb.
$
\end{itemize}
Let $L_Y(\tau)$ be the fundamental solution
of the quantum differential equation,
that is, the $\End(H^\bullet_\amb(Y; \bC))$-valued functions
satisfying
$
 \nabla_i^A L_Y(\tau) = 0
$
for $i = 1, \ldots, r$ and
$
 L_Y(\tau) = \id + O(\tau).
$
Since $Y$ is a K3 surface,
the quantum cup product $\circ_\tau$ coincides
with the ordinary cup product, and
the fundamental solution is given by
$
 L(\tau) = \exp(- \tau) \cup(-).
$
Let $\HAC = \Ker \nabla^A$ be the $\bC$-local system
associated with $\nabla^A$ and define the integral local subsystem
$\HAZ^\amb \subset \HAC$ as
\begin{align}
 \HAZ^\amb
  = \lc (2 \pi \sqrt{-1})^{-2} L_Y \lb \Gammahat_Y
   \cup (2 \pi \sqrt{-1})^{\frac{\deg}{2}} \ch(\iota^*\cE) \rb
      \relmid \scE \in K(X) \rc,
\end{align}
where the $\Gammahat_Y$ is defined in terms
of the Chern roots $\delta_1, \delta_2$ of the tangent bundle $TY$
as
$
 \Gammahat_Y := \Gamma(1+\delta_1) \Gamma(1+\delta_2)
$
\cite[Definition 6.3]{Iritani_QCP}.

The basis $\{ p_i \}_{i=1}^r$ of $\Pic X$ determines
a coordinate $q=(q_1,\cdots,q_r)$
on $\cM \cong \Pic X \otimes_\bZ \bCx$.
Let $u_i \in H^2(X; \bZ)$ be the Poincar\'{e} dual
of the toric divisor
corresponding to the one-dimensional cone
$\bR \cdot b_i \in \Sigma$ and
$v = u_1 + \cdots + u_m$
be the anticanonical class.
Givental's {\em $I$-function} is defined as the series
$$
 I_{X, Y}(q, z) = e^{p \log q / z}
  \sum_{d \in \Eff(X)} q^d \,
  \frac{
   \prod_{k=-\infty}^{\la d, v \ra} (v + k z)
   \prod_{j=1}^m \prod_{k=-\infty}^0
    (u_j + k z)}
  {\prod_{k=-\infty}^0
    (v + k z)
   \prod_{j=1}^m \prod_{k=-\infty}^{\la d, u_j \ra}
    (u_j + k z)},
$$
which gives a multi-valued map from an open subset of $\cM \times \bCx$
to the classical cohomology ring $H^\bullet(X; \bC[z^{-1}])$.
Givental's {\em $J$-function} is defined by
$$
 J_Y(\tau, z) = L_Y(\tau, z)^{-1}(1)
  = \exp(\tau/z).
$$
If we write
$$
 I_{X, Y}(q, z) = F(q) + \frac{G(q)}{z} + \frac{H(q)}{z^2} + O(z^{-3}),
$$
then Givental's mirror theorem
\cite{Givental_EGWI,
Givental_MTTCI,
Coates-Givental}
states that
\begin{align} \label{eq:Givental}
 \Euler(\omega_X^{-1}) \cup I_{X, Y}(q, z)
  = F(q) \cdot \iota_* J_Y(\varsigma(q), z),
\end{align}
where $\Euler(\omega_X^{-1}) \in H^2(X; \bZ)$
is the Euler class of the anticanonical bundle of $X$,
and the {\em mirror map}
$
 \varsigma(q) : \cM \to H^2_\amb(Y; \bC)
$
is a multi-valued map defined by
\begin{align}
 \varsigma(q) = \iota^* \lb \frac{G(q)}{F(q)} \rb.
\end{align}
The functions $F(q)$, $G(q)$ and $H(q)$ satisfy
the Gelfand--Kapranov--Zelevinsky hypergeometric differential equations,
and give periods for the B-model VHS
$
 (\HB, \nabla^B, \scrF_B^\bullet, Q_B).
$
Iritani \cite[Theorem 6.9]{Iritani_QCP} lifted
the mirror theorem \eqref{eq:Givental}
to an isomorphism
of integral variations of pure and polarized Hodge structures.

An important step
in the proof of \cite[Theorem 6.9]{Iritani_QCP}
is the identification,
given in the proof of \cite[Theorem 5.7]{Iritani_QCP},
of the monodromy of the B-model VHS
along an element $\ell$ of
$\pi_1(\bL_{\bCx}^\vee) \cong \bL^\vee$
in the neighborhood of the large complex structure limit point,
with the isometry
\begin{align} \label{eq:lb_tensor}
 (-) \otimes \iota^*(\cL^\vee) \colon \cN(\Yv_\alpha) \to \cN(\Yv_\alpha)
\end{align}
induced by the tensor product of the restriction to $\Yv_\alpha$
of the dual of the line bundle $\cL$ on $\Xv$
with $c_1(\cL) = \ell \in \bL^\vee \cong \Pic \Xv$.
Note that the isometry \eqref{eq:lb_tensor} lifts
to an autoequivalence of $D^b \coh \Yv$.
The relation between monodromy of period map and
autoequivalence of the derived category
goes back to \cite{Kontsevich_ENS98,Horja_HFMSTV}.

Givental's mirror theorem \eqref{eq:Givental}
(and its integral lift by Iritani) gives
an isomorphism of the A-model VHS and the B-model VHS
only in a neighborhood of the large radius limit point.
A global study of the period map
from the point of view of mirror symmetry
has been done
for the quartic mirror family of K3 surfaces
in \cite{MR3062592}
based on earlier works 
\cite{Nagura-Sugiyama,Narumiya-Shiga}.
Similar analysis for a couple of 2-parameter cases
has been performed in \cite{1403.5818}
based on earlier works
\cite{MR2962398,MR3160602}.

\section{Dolgachev conjecture}
 \label{sc:D_conj}

Let $(\Delta, \Deltav)$ be a polar dual pair
of three-dimensional reflexive polytopes.
Let further $X$ be a smooth toric weak Fano 3-fold
whose fan polytope is $\Deltav$, and
$\Xv$ be another smooth toric weak Fano 3-fold
whose fan polytope is $\Delta$.
In other words,
$\Delta$ is the Newton polytope
of $H^0(\scO_X(- K_X))$, and
$\Deltav$ is the Newton polytope
of $H^0(\scO_\Xv(- K_\Xv))$.
Let $Y \subset X$ and $\Yv \subset \Xv$ be
smooth anticanonical hypersurfaces.
Define $M_\Delta \subset H^2(Y;\bZ)$
as the primitive sublattice
generated by $H_{A, \, \bZ}^\amb(Y) \cap H^2(Y, \bZ)$,
and similarly for $M_\Deltav \subset H^2(\Yv, \bZ)$.

\begin{conjecture}[{Dolgachev
\cite[Conjecture (8.6)]{Dolgachev_MSK3}}]
 \label{cj:Dolgachev}
\ 
\begin{enumerate}
\item
There exist a lattice $\Mv_\Delta$
and an orthogonal decomposition
$M_\Delta^\bot = U \bot \Mv_\Delta$.
\item
There exists a primitive embedding
$
 M_\Deltav \subset \Mv_\Delta.
$
\item
The equality
$
 M_\Deltav = \Mv_\Delta
$
holds if and only if
$
 M_\Delta \cong \Pic Y.
$
\end{enumerate}
\end{conjecture}

Let $\HAZ(Y)$ and $\HBZ(Y)$ be
the primitive sublattices of $H^*(Y; \bZ)$
generated by $\HAZ^\amb(Y)$ and $\HBZ^\vc(Y)$
respectively.
Let further $M_\Delta^0$ be the orthogonal complement
of $M_\Delta$ inside $\NS(Y)$.
One has
\begin{align}
 \HAZ(Y) &= U \oplus M_\Delta
\end{align}
where $U$ is the unimodular hyperbolic plane
generated by the classes $[\cO_Y]$
and $[\cO_p]$
of the structure sheaf and a skyscraper sheaf.
One has sublattices
\begin{align}
\begin{split}
 \HAZ^\amb(Y) \bot M_\Delta^0 \bot \HBZ^\vc(Y)
  &\subset H^*(Y; \bZ), \\
 \HBZ^\vc(\Yv) \bot M_\Deltav^0 \bot \HAZ^\amb(\Yv)
  &\subset H^*(\Yv; \bZ)
\end{split}
\end{align}
and isomorphisms
\begin{align}
\begin{split}
 \Mir_Y^\bZ \colon \HAZ^\amb(Y) \simto \HBZ^\vc(\Yv), \\
 \Mir_\Yv^\bZ \colon \HAZ^\amb(\Yv) \simto \HBZ^\vc(Y).
\end{split}
 \label{eq:Mir1}
\end{align}
The ranks of the sublattices are given by
\begin{align}
 \rank \HAZ^\amb(Y) &= \#((\Delta^\circ)^{(1)} \cap \bsN) - 3, \\
 \rank \HBZ^\vc(Y) &= \#(\Delta^{(1)} \cap \bsM) - 3, \\
 \rank M_\Delta^0 &= \sum_{\gamma \in \Delta^{(1)}}
  \ell^*(\gamma) \ell^*(\gamma^*),
\end{align}
where $\Delta^{(1)}$ denotes the 1-skeleton
of the polytope $\Delta$, and
$\ell^*(\gamma)$ denotes the number of interior lattice points
of an interval $\gamma$
(\cite[\S 5.11]{Danilov-Khovansky},
\cite[Proposition 4.3.3]{Kobayashi_DW}).
\pref{cj:Dolgachev}.1 holds
since the isomorphisms \eqref{eq:Mir1}
lifts to the hyperbolic plane
$
 U = \bZ [\cO_\Yv] \oplus \bZ [\cO_p]
  \subset \HBZ(\Yv)
$
by \cite[Theorem 6.10]{Iritani_QCP}.
Conjectures \ref{cj:Dolgachev}.2 and
\ref{cj:Dolgachev}.3 hold
if the isomorphisms \eqref{eq:Mir1}
lifts to an isomorphism
\begin{align}
\begin{split}
 \HAZ(Y) &\simto \HBZ(\Yv), \\
 \HAZ(\Yv) &\simto \HBZ(Y)
\end{split}
\end{align}
of overlattices.

\section{Stability conditions}
 \label{sc:stab}

Stability conditions are introduced by Bridgeland
\cite{Bridgeland_SCKS}
motivated by stability of BPS D-branes
studied by string theorists
\cite{Douglas_DBHMSS}.

\begin{definition}[{\cite[Definition 1.1]{Bridgeland_SCTC}}]
A \emph{stability condition}
on a triangulated category $\scT$
is a pair $\sigma = (Z,\scP)$
consisting of
\begin{itemize}
 \item
a group homomorphism $Z \colon K(\scT)\to\bC$, and
 \item
full additive subcategories
$\scP(\phi)$ for each $\phi\in\bR$
\end{itemize}
satisfying the following conditions:
\begin{enumerate}
 \item
If $0\ne E\in \scP(\phi)$, then $Z(E)=m(E)\exp(i\pi\phi)$ for some
$m(E)\in\bR_{>0}$,
\item
for all $\phi\in\bR$, $\scP(\phi+1)=\scP(\phi)[1]$,
\item
for $A_j\in\scP(\phi_j)$ $(j=1,2)$ with $\phi_1>\phi_2$,
one has
$\Hom_{\scT}(A_1,A_2)=0$,
\item
for every non-zero object $E\in\scT$, there is a finite sequence of
real numbers
\[\phi_1>\phi_2>\cdots>\phi_n\]
and a collection of triangles
\begin{equation} \label{eq:HN_filtration}
\psmatrix[colsep=.5, rowsep=.7]
  0 = \Rnode{a}{E_0} & & E_1 & & E_2 & \cdots & E_{n-1} & &  
\Rnode{n}{E_n}=E \\
        & A_1 & & A_2 &   &        &       & A_n & \\
\endpsmatrix
\psset{nodesep=3pt,arrows=->}
\ncline{a}{1,3}
\ncline{1,3}{1,5}
\ncline{1,6}{1,7}
\ncline{1,7}{n}
\ncline{1,3}{2,2}
\ncline{1,5}{2,4}
\ncline{n}{2,8}
\psset{arrows=-}
\ncline{1,5}{1,6}
\psset{arrows=->,linestyle=dashed}
\ncline{2,2}{a}
\ncline{2,4}{1,3}
\ncline{2,8}{1,7}
\end{equation}
with $A_j\in\scP(\phi_j)$ for all $j$.
\end{enumerate}
\end{definition}

The homomorphism $Z$ is called the {\em central charge},
and the diagram \eqref{eq:HN_filtration}
is called the {\em Harder-Narasimhan filtration}.
It follows from the definition
that $\scP(\phi)$ is an abelian category,
and its non-zero object $\scE \in \scP(\phi)$ is said to be
{\em semistable of phase $\phi$}.
An object $\scE$ is said to be {\em stable}
if it is a simple object of $\scP(\phi)$,
i.e., there are no proper subobjects
of $\scE$
in $P(\phi)$.
By \cite[Proposition 5.3]{Bridgeland_SCTC},
giving a stability condition
on a triangulated category $\scT$
is equivalent to
giving a bounded $t$-structure on $\scT$
and a {\em stability function}
on its heart
with the {\em Harder-Narasimhan property}.
For the definitions of
a stability function
and the Harder-Narasimhan property,
see \cite[\S 2]{Bridgeland_SCTC}.

Assume that the numerical Grothendieck group
$\cN(\cN)$ of finite rank,
and take a norm $\| \cdot \|$
on the finite-dimensional vector space
$\cN(\cT) \otimes \bR$.
A stability condition is \emph{numerical}
if the central charge $Z \colon K(\cT) \to \bC$
factors through the numerical Grothendieck group.
A numerical stability condition is said to satisfy
the \emph{support condition}
\cite{Kontsevich-Soibelman_SSMDTICT}
if there is a positive constant $K$
such that for any $E \in \cP(\phi)$,
one has
\begin{align}
 |Z(E)| \ge K \| E \|.
\end{align}
The set of numerical stability conditions
satisfying the support condition
is denoted by $\Stab \scT$.
There is a natural topology
on $\Stab\scT$
such that the forgetful map
\begin{align} \label{eq:cZ}
\begin{array}{rccc}
 \cZ : & \Stab \scT & \to & \Hom(\cN(\cT), \bC) \\
 & \rotatebox{90}{$\in$} & & \rotatebox{90}{$\in$} \\
 & (Z, \scP) & \mapsto & Z
\end{array}
\end{align}
is a local homeomorphism.
This local homeomorphism induces a structure
of a complex manifold on $\Stab \cT$.

Since the definition of $\Stab \scT$
depends only on the triangulated structure of $\scT$,
the group $\Aut \scT$
of autoequivalences of $\scT$ as a triangulated category
acts naturally on $\Stab \scT$ from the left;
for $\sigma = (Z, \scP) \in \Stab \scT$
and
$\Phi \in \Aut \scT$,
\begin{align}
 \Phi(\sigma) = (Z \circ \Phi_*^{-1}, \Phi(\scP))
\end{align}
where
$
 \Phi_* : \cN(\scT) \to \cN(\scT)
$
is the automorphism
induced by $\Phi$.
This action commutes with the right action
of the universal cover
$\GLt(2, \bR)$
of the general linear group
$GL^+(2, \bR)$
with positive determinant,
which ``rotates'' the central charge
\cite[Lemma 8.2]{Bridgeland_SCTC}.

If $\cT = D(Y)$ is the bounded derived category
of coherent sheaves on a smooth projective variety $Y$,
then $\Stab \cT$ will be denoted by $\Stab(Y)$.
A stability condition is \emph{geometric}
if all skyscraper sheaves $\cO_x$ are stable
of the same phase.
When $Y$ is a K3 surface,
then there is a distinguished connected component
$\Stab^\dagger(Y)$ containing geometric stability conditions.
Let $\Stab^*(Y)$ be the union of connected components
of $\Stab(Y)$
which are images of $\Stab^\dagger(Y)$
by the action of $\Aut D(Y)$.

The derived categories $D(Y)$ and $D(Y')$
of K3 surfaces $Y$ and $Y'$ are equivalent
as a triangulated category
if and only if there exists a Hodge isometry
between the transcendental lattices
of $Y$ and $Y'$
\cite[Theorem 3.3]{Orlov_EDCKS}.
An autoequivalence $\Phi \in \Aut D(Y)$
induces a Hodge isometry $\Phi_* : H^*(Y; \bZ) \to H^*(Y; \bZ)$
making the following diagram commutative;
\begin{align}
\begin{CD}
 D^b \coh Y @>{\Phi}>> D^b \coh Y \\
  @V{\ch(-) \sqrt{\td_Y}}VV @V{\ch(-) \sqrt{\td_Y}}VV \\
 H^*(Y; \bZ) @>{\Phi_*}>> H^*(Y; \bZ).
\end{CD}
\end{align}
This induces a group homomorphism
\begin{align}
 (-)_* \colon \Aut D(Y) \to \Aut H^*(Y; \bZ)
\end{align}
to the group of Hodge isometries,
%
%
%
whose image is the index 2 subgroup
\begin{align}
 \Aut^+ H^*(Y; \bZ) \subset \Aut H^*(Y;\bZ)
\end{align}
of Hodge isometries
preserving the orientation of positive-definite 4-planes
\cite{MR2047679, MR2553878}.
The kernel of this homomorphism
will be denoted by $\Aut^0 D(Y)$.

Consider the space
\begin{align}
 \cP = \{ \Omega \in \cN(Y) \otimes \bC
  \mid (\Omega, \Omegabar) > 0 \},
\end{align}
which can be identified with the space of pairs
$(\Re \Omega, \Im \Omega)$ of elements in $\cN(Y) \otimes \bR$
which span a positive-definite 2-planes in $\cN(Y) \otimes \bR$.
The space $\cP$ consists of two connected components.
We write the connected component
containing the vector
$(1, \sqrt{-1} \omega, - \frac{1}{2} \omega^2)$
for an ample class $\omega \in \NS(Y) \otimes \bR$
as $\cP^+(Y)$.
The complement of the union of reflection hyperplanes
for all $(-2)$-elements is denoted by
\begin{align}
 \cP_0^+(Y) = \cP^+(Y) \setminus \bigcup_{\delta \in \Delta(Y)} \delta^\bot.
\end{align}

\begin{theorem}[{\cite[Theorem 1.1]{Bridgeland_SCK3}}]
The map \eqref{eq:cZ} induces a Galois covering
\begin{align}
 \cZ \colon \Stab^\dagger(Y) \to \cP_0^+(Y)
\end{align}
whose group of deck transformations
is the subgroup of $\Aut^0 D(Y)$
which preserves the connected component $\Stab^\dagger(Y)$.
\end{theorem}

\pref{cj:Bridgeland} below implies that
$
 \Aut^0 D(Y) \cong \pi_1(\cP_0^+(Y)).
$

\begin{conjecture} \label{cj:Bridgeland}
$\Stab^*(Y)$ is connected and simply connected.
\end{conjecture}

Let
\begin{align}
 \cQ(Y) = \{ \Omega \in \cP(Y) \mid (\Omega, \Omega) = 0 \}
\end{align}
be a quadric hypersurface of $\cP(Y)$
and set $\cQ_0^+(Y) = \cQ(Y) \cap \cP_0^+(Y)$.
The group $\GL_2^+(\bR)$ acts freely on $\cP_0^+(Y)$
by rotating $(\Re \Omega, \Im \Omega)$.
The group $\bCx$ acts as a subgroup of $\GL_2^+(\bR)$
by rescaling $\Omega$.
Since $\GL_2^+(\bR)$-orbits in $\cP_0^+(Y)$ intersects $\cQ_0^+(Y)$
in a unique $\bCx$-orbit,
one has a homeomorphism
\begin{align}
 \cP_0^+(Y) / \GL_2^+(\bR) \cong \cQ_0^+(Y) / \bCx,
\end{align}
which induces a group isomorphism
\begin{align}
 \pi_1(\cP_0^+(Y) / \GL_2^+(\bR))
  \cong \pi_1(\cQ_0^+(Y) / \bCx).
\end{align}
Since
$
 \pi_1(\GL_2^+(\bR))
  \cong \pi(\bCx)
  \cong \pi_1(S^1)
  \cong \bZ,
$
this yields a group isomorphism
\begin{align} \label{eq:PQ}
 \pi_1(\cP_0^+(Y)) \cong \pi_1(\cQ_0^+(Y)).
\end{align}
Let
\begin{align} \label{eq:AutCYH}
 \Aut_\CY^+ H^*(Y;\bZ) \subset \Aut^+ H^*(Y;\bZ)
\end{align}
be the subgroup consisting of Hodge isometries
which acts trivially on $H^{2,0}(Y) \subset H^*(Y;\bZ) \otimes \bC$.
Such action is usually called \emph{symplectic},
where CY is the notation in \cite{Bayer-Bridgeland}
standing for \emph{Calabi-Yau}.
An element $\phi \in \Aut^+ H^*(Y;\bZ)$ lies in $\Aut_\CY^+ H^*(Y;\bZ)$
if and only if $\phi$ acts trivially
on the transcendental lattice $T(Y)$:
The `if' part is clear since $H^{2,0}(Y) \subset T(Y) \otimes \bC$,
and the `only if' part follows
from the fact that for any $\alpha \in T(Y)$,
the element $\phi(\alpha) - \alpha$ is integral and
orthogonal to both $H^{2,0}(Y)$ and $\cN(Y)$, and hence zero.
It follows that $\Aut_\CY^+ H^*(Y;\bZ)$ is isomorphic
to the index 2 subgroup of the group
$\Aut \cN(Y)$ of isometries of $\cN(Y)$
preserving orientations of positive-definite 2-planes.

Define $\Aut_\CY(D(Y))$ by
\begin{align} \label{eq:AutCYD}
 \Aut_\CY(D(Y))
  = \{ \Phi \in \Aut(D(Y)) \mid
   \Phi_* \in \Aut_\CY^+ H^*(Y;\bZ) \}.
\end{align}
\pref{cj:Bridgeland} gives 
an exact sequence
\begin{align}
 1
  \to \pi_1(\cP_0^+(Y))
  \to \Aut D(Y)
  \to \Aut^+ H^*(Y;\bZ)
  \to 1,
\end{align}
which together with \eqref{eq:PQ},
\eqref{eq:AutCYH} and \eqref{eq:AutCYD}
induces a short exact sequence
\begin{align}
 1
  \to \pi_1(\cQ_0(Y))
  \to \Aut_\CY D(Y)
  \to \Aut_\CY^+ H^*(Y;\bZ)
  \to 1
\end{align}
and an isomorphism
\begin{align}
 \pi_1^\orb([\cQ_0(Y)/\Aut_\CY^+ H^*(Y;\bZ)])
  \cong \Aut_\CY D(Y).
\end{align}
Since the natural action of $\bCx$
on $\cQ_0^+(Y)$ is free and
commutes with the action of $\Aut_\CY^+ H^*(Y)$,
the orbifold quotient $[\cQ_0(Y)/\Aut_\CY^+ H^*(Y;\bZ)]$
is a principal $\bCx$-bundle
over the moduli space
$\cM_0(Y) = [\cD_0^+(Y)/ \Aut_\CY^+ H^*(Y;\bZ)]$
of ample $\NS(Y)$-polarized K3 surfaces.
The associated long exact sequence of homotopy groups
\begin{align}
 1 \to \pi_1(\bCx)
  \to \pi_1^\orb([\cQ_0(Y)/\Aut_\CY^+ H^*(Y;\bZ)])
  \to \pi_1^\orb(\cM_0(Y))
  \to 1
\end{align}
gives
\begin{align}
 1
  \to \bZ
  \to \Aut_\CY D(Y)
  \to \pi_1^\orb(\cM_0(Y))
  \to 1.
\end{align}
The map
$
 \bZ \to \Aut_\CY D(Y)
$
sends $1 \in \bZ$
to the shift $[2] \in \Aut_\CY D(Y)$,
so that \pref{cj:Bridgeland} implies
(and is in fact equivalent to)
the isomorphism
\begin{align}
 \pi_1^\orb(\cM_0(Y)) \cong \Aut_\CY D(Y)/[2].
\end{align}

\section{Borcea-Voisin mirrors}
 \label{sc:BV}


An involution $\iota \colon S \to S$ on a K3 surface $S$
is said to be \emph{symplectic}
if the induced map $\iota^* \colon H^{2,0}(S) \to H^{2,0}(S)$
on $H^{2,0}(S)$ is the identity map,
and \emph{anti-symplectic} if it is minus the identity map.
A \emph{2-elementary K3 surface} is a pair $(K, \iota)$
of a K3 surface and an anti-symplectic involution $\iota$.
Any 2-elementary K3 surface is algebraic.

Decompose $H^2(S; \bC)$ into the $(+1)$-and $(-1)$-eigenspaces
$H^2(S; \bC)^+$ and $H^2(S; \bC)^-$
of the action $\iota^* \colon H^2(S; \bC) \to H^2(S; \bC)$,
and set $H^2(S; \bZ)^\pm = H^2(S; \bC)^\pm \cap H^2(S; \bZ)$.
The signatures of $H^2(S; \bZ)^+$ and $H^2(S; \bZ)^-$ are given by
$(1, r-1)$ and $(2, 20-r)$ respectively,
where $r$ is the rank of $H^2(S; \bZ)^+$.

A lattice $M$ of rank $r$ is \emph{2-elementary}
if $A_M = M^\vee/M \cong (\bZ/2\bZ)^a$,
where $a$ is the minimal number of generators of $A_M$.
The induced quadratic form $q_M : A_M \to \bQ/2 \bZ$
is called the \emph{discriminant form}.
The \emph{parity} of $q_M$ is defined by
\begin{align}
 \delta =
\begin{cases}
 0 & q_M(A_M) \subset \bZ, \\
 1 & \text{otherwise}.
\end{cases}
\end{align}
The triple $(r, a, \delta)$ is called the \emph{main invariant} of $M$.

\begin{theorem}[{\cite[Theorem 4.3.2]{MR556762}}]
The isometry class of a hyperbolic even 2-elementary lattice is determined
by the main invariant.
\end{theorem}

If $(S, \iota)$ is a 2-elementary K3 surface,
then $H^2(S; \bZ)^+$ s a 2-elementary lattice.
The \emph{main invariant} of a 2-elementary K3 surface $(S, \iota)$
is defined as the main invariant
of the 2-elementary lattice $H^2(S; \bZ)^+$.

\begin{theorem}[{\cite[Theorem 4.2.2]{MR556762}}]
Let $(S, \iota)$ be a 2-elementary K3 surface
with main invariant $(r,a,\delta)$.
Then the fixed locus $S^\iota \subset S$ is given as follows:
\begin{itemize}
 \item
If $(r, a, \delta) = (10, 10, 0)$,
then $S^\iota = \emptyset$. 
 \item
If $(r,a,\delta)=(10,8,0)$,
then $S^\iota$ is the union of two elliptic curves.
 \item
In all other cases,
$S^\iota$ is the disjoint union
$S^\iota = C^g \coprod E_1 \coprod \cdots \coprod E_k$
of a curve $C^g$ of genus $g$
and rational curves $E_1, \ldots, E_k$
with
\begin{align}
 g = 11 - \frac{r+a}{2}, \quad
 k = \frac{r-a}{2}.
\end{align}
\end{itemize}
\end{theorem}

Let $\iota \colon E \to E$ be an involution on an elliptic curve $E$.
The induced map $\iota^* \colon H^{1,0}(E) \to H^{1,0}(E)$
is either $\id$ or $- \id$.
The quotient $E/ \la \iota \ra$ is isomorphic
to an ellitptic curve if $\iota^* = \id$, and
to $\bP^1$ if $\iota^* = - \id$.

Let $(S, \iota_1)$ be a 2-elementary K3 surface,
and $(E, \iota_2)$ be an elliptic curve with an involution
with $E / \la \iota_2 \ra \cong \bP^1$.
Then the product $\iota = \iota_1 \times \iota_2$
is an involution on $S \times E$
such that the induced map
$
 \iota^* \colon H^{3,0}(S \times E) \to H^{3,0}(S \times E)
$
is the identity map.
If we write the irreducible components
of the fixed loci of $\iota_1$ and $\iota_2$
as $C_1, \ldots, C_n$ and $p_1, \ldots p_4$ respectively,
then the fixed point of $\iota$
consists of $C_i \times p_j$ for $i=1, \ldots, n=k+1$ and $j=1,\ldots,4$.
The involution $\iota$ lifts naturally
to an involution $\iotatilde$
on the blow-up $(S \times E)^\sim$
of $S \times E$ along these $4 n$ curves,
and the quotient will be denoted by $Y=(S \times E)^\sim/\la \iotatilde \ra$.

\begin{theorem}[\cite{MR1416355,MR1265318}]
$Y$ is a smooth Calabi-Yau 3-fold with
\begin{align}
 h^{1,1}(Y) = 11+5n-g, \quad
 h^{2,1}(Y)=11+5g-n.
\end{align}
\end{theorem}

The classification of hyperbolic
even 2-elementary lattices
admitting primitive embedding in the K3 lattice
is given in \cite[Table 1]{MR556762}.
One can immediately see from the table
that if one restricts to $g>0$ and $(r,a,\delta) \ne (14,6,0)$,
then there is a 2-elementary K3 surface $(S, \iota)$
with main invariants $(r,a,\delta)$
if and only if there is a 2-elementary K3 surface $(\Sv, \iotav)$
with main invariants $(\rv, \av, \deltav) = (20-r,a,\delta)$.
Take an elliptic curve $(\Ev, \iotav_2)$
with an involution such that $\Ev/\iotav_2 \cong \bP^1$.
The Calabi-Yau 3-fold
associated with $((\Sv, \iotav_1), (\Ev, \iotav_2))$
will be denoted by
$\Yv = (\Sv \times \Ev)^\sim/\la (\iotav_1, \iotav_2) \ra$.
Then one has
\begin{align}
 h^{1,1}(Y) = h^{2,1}(\Yv), \quad
 h^{2,1}(Y) = h^{1,1}(\Yv).
\end{align}
The pair $(Y, \Yv)$ is called
a \emph{Borcea-Voisin mirror pair}.

\bibliographystyle{amsalpha}
\bibliography{bibs}

\noindent
Kazushi Ueda

Department of Mathematics,
Graduate School of Science,
Osaka University,
Machikaneyama 1-1,
Toyonaka,
Osaka,
560-0043,
Japan.

{\em e-mail address}\ : \  kazushi@math.sci.osaka-u.ac.jp
\ \vspace{0mm} \\

\end{document}